\documentclass [11pt,a4paper]{article}
\usepackage[ansinew]{inputenc}
\usepackage{amsthm,amsmath,amssymb}
\usepackage{graphicx}
\usepackage{cite}
\usepackage{mathrsfs}

\usepackage{color}

\newtheorem{theorem}{\textbf{Theorem}}[section]
\newtheorem{lemma}{\textbf{Lemma}}[section]
\newtheorem{proposition}{\textbf{Proposition}}[section]
\newtheorem{corollary}{\textbf{Corollary}}[section]
\newtheorem{remark}{\textbf{Remark}}[section]
\newtheorem{definition}{\textbf{Definition}}[section]

\allowdisplaybreaks[4]

\def\be{\begin{equation}}
\def\ee{\end{equation}}
\def\bea{\begin{eqnarray}}
\def\eea{\end{eqnarray}}
\def\bt{\begin{theorem}}
\def\et{\end{theorem}}
\def\bl{\begin{lemma}}
\def\el{\end{lemma}}
\def\br{\begin{remark}}
\def\er{\end{remark}}
\def\bp{\begin{proposition}}
\def\ep{\end{proposition}}
\def\bc{\begin{corollary}}
\def\ec{\end{corollary}}
\def\bd{\begin{definition}}
\def\ed{\end{definition}}

\def\p{\partial}

\def\non{\nonumber }
\def\wha{\widehat}
\def\wt{\widetilde}

\definecolor{coloras}{rgb}{0.,0.67,0}
\begin{document}

\title{Long-time behavior \\ for a hydrodynamic model on nematic liquid
crystal flows \\ with asymptotic stabilizing boundary condition and external force}

\author{
{\sc Maurizio Grasselli}\ \footnote{Dipartimento di Matematica,
Politecnico di Milano, Milano 20133, Italy, {\it
maurizio.grasselli@polimi.it}}\ \
and\ {\sc Hao Wu}\footnote{Corresponding author. School of Mathematical
Sciences and Shanghai Key Laboratory for Contemporary Applied Mathematics,
Fudan University, Shanghai 200433, China, \textit{haowufd@yahoo.com}}
}

\date{\today}

\maketitle


\begin{abstract}\noindent
In this paper, we consider a simplified Ericksen--Leslie
model for the nematic liquid crystal flow. The evolution system consists of the Navier--Stokes
equations coupled with a convective Ginzburg--Landau type equation
for the averaged molecular orientation. We
suppose that the Navier--Stokes
equations are characterized by a no-slip boundary condition and a
time-dependent external force $\mathbf{g}(t)$, while the equation for the molecular director is subject to a time-dependent Dirichlet boundary condition $\mathbf{h}(t)$.
 We show that, in 2D,
each global weak solution converges to a single stationary state
when $\mathbf{h}(t)$ and $\mathbf{g}(t)$ converge to a
time-independent boundary datum $\mathbf{h}_\infty$ and
$\mathbf{0}$, respectively. Estimates on the
convergence rate are also obtained. In the 3D case, we prove that global weak solutions
are eventually strong so that results similar to the 2D case can be proven.
We also show the existence of global strong solutions, provided that either the
viscosity is large enough or the initial datum is close to a given
equilibrium.\medskip


\end{abstract}

\section{Introduction}
\label{Intro}
We consider the following hydrodynamical model for the
flow of nematic liquid crystals
 \bea
 \mathbf{v}_t+\mathbf{v}\cdot\nabla \mathbf{v}-\nu \Delta \mathbf{v}
 +\nabla \pi&=&-\lambda
 \nabla\cdot(\nabla \mathbf{d}\odot\nabla \mathbf{d})+\mathbf{g}(t),\label{1}\\
 \nabla \cdot \mathbf{v} &=& 0,\label{2}\\
 \mathbf{d}_t+\mathbf{v}\cdot\nabla \mathbf{d}&=&\eta(\Delta \mathbf{d}-\mathbf{f}(\mathbf{d})),\label{3}
 \eea
in $\Omega \times\mathbb{R}^+$, where $\Omega \subset \mathbb{R}^n\
(n=2,3)$ is a bounded domain with sufficiently smooth boundary
$\Gamma$, $\mathbf{v}=(v_1,...,v_n)^{tr}$ is the velocity field of
the flow and $\mathbf{d}=(d_1,...,d_n)^{tr}$ represents the averaged
macroscopic/continuum molecular orientations in $\mathbb{R}^n \
(n=2,3)$. $\pi$ is a scalar function representing the pressure
(including both the hydrostatic and the induced elastic part from
the orientation field). The external volume force is represented by
$\mathbf{g}$. The positive constants $\nu, \lambda$ and $\eta$ stand
for viscosity, the competition between kinetic energy and potential
energy, and macroscopic elastic relaxation time (Deborah number) for
the molecular orientation field. $\nabla \mathbf{d}\odot \nabla
\mathbf{d}$ denotes the $n\times n$ matrix whose $(i,j)$-th entry is
given by $\nabla_i \mathbf{d}\cdot \nabla_j \mathbf{d}$, for $1\leq
i,j\leq n$. We assume that
$\mathbf{f}(\mathbf{d})=\nabla_{\mathbf{d}} F(\mathbf{d})$ for some
smooth bounded function $F:\mathbb{R}^n\rightarrow \mathbb{R}$. In
particular, one uses the Ginzburg--Landau approximation
$\mathbf{f}(\mathbf{d})=\frac{1}{\epsilon^2}(|\mathbf{d}|^2-1)\mathbf{d}$
to relax the nonlinear constraint $|\mathbf{d}|=1$ on molecule length (cf.
\cite{LL95,LLW}).

System \eqref{1}--\eqref{3} was firstly proposed in \cite{lin1} as a
simplified approximate system of the original Ericksen--Leslie model
for the nematic liquid crystal flows (cf. \cite{E1,Le}).
Well-posedness of the autonomous version of system
\eqref{1}--\eqref{3} (namely, with $\mathbf{g}= \mathbf{0}$, no-slip
boundary condition for $\mathbf{v}$, and time-independent Dirichlet
boundary condition for $\mathbf{d}$) has been analyzed in
\cite{LL95} (see also \cite{FO09,LL96} and, for different
boundary conditions, \cite{LS01}). For numerical approximation
we refer to \cite{LL06,Liu2,Liu1}. Problem
\eqref{1}--\eqref{3} has also been investigated on a Riemannian
manifold in \cite{S01}, where the existence of a global attractor in
the 2D case was proven. As far as the long-time behavior of the
single trajectory is concerned, in \cite{LL95}, a natural question
on the uniqueness of asymptotic limit for global solutions (to the
autonomous system) was raised. This question was answered in \cite{W10}, where it is proven that each trajectory converges
to a single steady state (cf. \cite{PRS,wuxuliu} for some generalization). The proof is based on a suitable
{\L}ojasiewicz--Simon type inequality (see
\cite{S83}, cf. also \cite{Hu} and references cited therein).

 The technically more challenging case of time-dependent
Dirichlet boundary conditions for $\mathbf{d}$ has been recently
analyzed in \cite{B,C06,C09,G09}. For instance, under proper
assumptions on the time-dependent boundary condition and assuming
that $\mathbf{g}=\mathbf{0}$, the existence of global weak solution,
the existence of global regular solution for viscosity coefficient
big enough, and the weak/strong uniqueness were obtained in
\cite{C09}. Regularity criteria for solutions in the 3D case can be
found in \cite{G09}. Besides, the presence of a time-dependent
external force is allowed in \cite{B} and existence of global and
exponential attractors is proven in the 2D case. In this paper, we
want to extend the results of \cite{W10} to the non-autonomous case
treated in \cite{B}. Thus we consider system
\eqref{1}--\eqref{3} subject to the boundary conditions
 \be
 \mathbf{v}(x,t)=\mathbf{0},\quad \mathbf{d}(x,t)=\mathbf{h}(x,t),\qquad (x, t)\in \Gamma\times
 \mathbb{R}^+,
 \label{4}
 \ee
and the initial conditions
 \be
 \mathbf{v}|_{t=0}=\mathbf{v}_0(x) \ \ \text{with}\ \nabla\cdot \mathbf{v}_0=0,\quad
 \mathbf{d}|_{t=0}=\mathbf{d}_0(x),\qquad x\in \Omega.\label{5}
 \ee
In the 2D case, we prove that each weak/strong solution converges to a single
stationary state when $\mathbf{h}(t)$ and $\mathbf{g}(t)$ converge to a time-independent boundary datum $\mathbf{h}_\infty$
and $\mathbf{0}$, respectively. In the 3D case, we first show the eventual regularity of global weak solutions, and the existence of global strong solutions provided that either the
viscosity is large enough or the initial datum is close to a given
equilibrium. Then an analogous result on the long-time behavior as in 2D is also obtained. In both cases,
we provide an estimate on the convergence rate.

Before ending this section, we state some key ingredients of the
present paper. System \eqref{1}--\eqref{5} is non-autonomous due to
the time-dependent boundary data $\mathbf{h}$ and external force
$\mathbf{g}$. This brings some additional difficulties into our
subsequent proofs. First, in order to obtain the energy inequalities
that play crucial roles in the proof of well-posedness as well as in
the long-time behavior of global solutions (cf. Lemmas \ref{BEL},
\ref{H2D}, \ref{H3D}, \ref{H3Ds}), we have to introduce proper
lifting functions (cf. \eqref{LE} and \eqref{LP} below). The idea
was first used in \cite{C06, C09}, but the lifting functions
introduced in this paper are different from those in \cite{C09}.
This is due to the fact that we need some specific energy inequalities which not only yield uniform estimates of the
solutions, but also provide estimates of the convergence rate (cf.
Section 4). The second issue regards the application of the \L
ojasiewicz--Simon approach (cf. \cite{S83}) which has been shown to
be very useful in the study of long-time behavior of global
solutions to nonlinear evolution equations
(see, for instance, \cite{J981, HJ01, Hu,WGZ1,W10,Z04} and
references therein). In particular, convergent results related to
various evolution equations with asymptotically autonomous source
terms were established, e.g., in
\cite{CJ03,GPS1,HT,Hu}. However, our current case is much more
complicated than the previous cases, because the \L ojasiewicz--Simon
inequality involves the vector $\mathbf{d}$ that is subject to a
time-dependent boundary datum.
To overcome this difficulty, we derive an extended \L
ojasiewicz--Simon type inequality for vector functions with
arbitrary nonhomogeneous Dirichlet boundary data, which is
associated with the lifted energy (cf. Corollary \ref{ELS}). This
generalizes the results in \cite{Hu, W10} and should have its own
interest. Third, in the 3D case, we also apply the \L
ojasiewicz--Simon approach to prove the existence of global strong
solutions provided that the initial datum is close to a local
minimizer of the elastic energy and the non-autonomous terms are properly
small perturbations of their asymptotic limits (cf. Section 5). Then
we further discuss the stability of these energy minimizers. This extends
the previous results in \cite{LL95,W10} for the autonomous system,
where the initial datum was required to be sufficiently close to a
global energy minimizer. For the stability of the general Ericksen--Leslie system \cite{LL01}, we refer to the recent work \cite{WXL13}.

The remaining part of the paper is organized as follows. The next
section is devoted to report some existence and uniqueness results
and  basic \emph{a priori} estimates for the solutions. The extended
{\L}ojasiewicz--Simon inequality we need is derived in Section
\ref{Sec4}. In Section \ref{Sec5} we show the convergence of each
global weak/strong solution to a single steady state and provide
uniform estimates on the convergence rate in 2D. Results in 3D are
presented in Section \ref{Sec6}. In particular, we study the eventual regularity of global weak solutions as well as the
well-posedness when the initial data are close to local minimizers of
the elastic energy. Long-time convergence of global solutions and stability of such minimizers are also proved. In the final Section
\ref{App}, some useful
properties of the lifting functions are reported.

\section{Preliminaries: well-posedness and \emph{a priori} estimates}
\label{Sec2} \setcounter{equation}{0}
Without loss of generality, from now on we set $\lambda=\eta=1$.
Let us introduce the function spaces we shall work with.  As
usual, $L^p(\Omega)$ and $W^{k,p}(\Omega)$ stand for the Lebesgue
and the Sobolev spaces of real valued functions, with the convention
that $H^k(\Omega)= W^{k,2}(\Omega)$. The spaces of vector-valued
functions are denoted by bold letters, correspondingly. Without any
further specification, $\|\cdot\|$ stands for the norm in
$L^2(\Omega)$ or $\mathbf{L}^2(\Omega)$. This norm is induced by the
scalar inner product $(u,v)=\int_\Omega uv dx$, where for vector
valued functions the product $uv$ is replaced by the Euclidean inner
product $\textbf{u}\cdot \textbf{v}$. We set, as usual, \be
\mathbf{H}=\overline{\mathcal{V}}^{\mathbf{L}^2(\Omega)},\quad
\mathbf{V}=\overline {\mathcal{V}}^{\mathbf{H}_0^1(\Omega)},\quad
\text{where}\ \mathcal{V}=\left\{\mathbf{v}\in
C_0^\infty(\Omega,\mathbb{R}^n)
 :\, \nabla\cdot \mathbf{v}=0\right\}.\non
 \ee
 For any Banach space $B$, we denote its dual space by $B^*$. In
 particular, we denote the dual space of $\mathbf{H}_0^1(\Omega)$ by
 $\mathbf{H}^{-1}(\Omega)$.

In the following text, we will use the regularity result for Stokes
problem (see, e.g., \cite{Te1})
  \bl \label{S}
  For the Stokes operator $S: D(S)= \mathbf{ V}\cap\mathbf{H}^2(\Omega) \to
  \mathbf{ H} $ defined by
  $$S\mathbf{u}=-\Delta \mathbf{u} +\nabla \pi \in \ \mathbf{H}, \quad
  \forall \mathbf{u} \in D(S),$$
  it holds
  $$\|\mathbf{u}\|_{\mathbf{H}^2}+\|\pi\|_{H^1\setminus\mathbb{R} }
  \leq C\|S \mathbf{u}\|, \quad \forall \mathbf{u} \in D(S),$$
for some positive constant $C$ only depending on $\Omega$ and the spatial dimension.
 \el

We begin to report the existence of a weak solution (see
\cite[Corollary 1.1, Theorem 1.4]{B}).
 \bp \label{we} Suppose $n=2,3$. For any given $T>0$, assume
 \begin{align}
 \label{hyp1}
&\mathbf{g}\in L^2(0, T; \mathbf{V}^*),\\
\label{hyp2}
&\mathbf{h}\in L^2(0, T;\mathbf{H}^\frac32(\Gamma)),\\
\label{hyp3}
&\mathbf{h}_t\in L^2(0,T; \mathbf{H}^{-\frac12}(\Gamma))\\
\label{hyp4}
&|\mathbf{h}|_{\mathbb{R}^n}\leq 1, \quad \textrm{ a.e. on }
\Gamma\times [0, T],\\
\label{hyp5}
&\mathbf{d}_0|_\Gamma=\mathbf{h}|_{t=0}.
\end{align}
Then for any $(\mathbf{v}_0, \mathbf{d}_0)\in \mathbf{H}\times
 \mathbf{H}^1(\Omega)$ with $|\mathbf{d}_0|_{\mathbb{R}^n}\leq 1$
 almost everywhere in $\Omega$, problem \eqref{1}--\eqref{5} admits a
 weak solution $(\mathbf{v}, \mathbf{d})$  such that
\begin{align}
&\mathbf{v}\in L^\infty(0, T; \mathbf{H})\cap L^2(0, T; \mathbf{V}),\non\\
&\mathbf{d}\in L^\infty(0, T; \mathbf{H}^1(\Omega)) \cap L^2(0, T; \mathbf{H}^2(\Omega)),\non\\
\label{max} &|\mathbf{d}(x,t)|_{\mathbb{R}^n}\leq 1, \quad { \rm a.e.\
on }\ \Omega\times [0, T].
\end{align}
If $n=2$, then the weak solution $(\mathbf{v}, \mathbf{d})$ to problem
\eqref{1}--\eqref{5} is unique. Moreover, we have $ (\mathbf{v}, \mathbf{d})\in C([0, T];
\mathbf{H}\times \mathbf{H}^1(\Omega))$.
 \ep

\br
The weak maximum principle \eqref{max} plays an important role in the analysis of system \eqref{1}--\eqref{5} (cf. \cite{LL95} for the autonomous case).
We recall that system \eqref{1}--\eqref{5} is a simplified version of the Ericksen--Leslie system for the liquid crystal flow of nematic type, in which the molecule is assumed to be ``small" such that the stretching and rotating effects in the fluid are neglected. In particular, when the stretching effect is taken into account (cf. \cite{sunliu}), the weak maximum principle \eqref{max} fails. The lack of control of the $L^\infty_t L^\infty_x$-norm of $\mathbf{d}$ brings extra difficulties in the analysis. For instance, in this case, it is not clear how to define weak solutions (compare with \cite{LL95,LL01}). We refer to \cite{CR, sunliu, GrasWu, PRS, wuxuliu} for extensive studies (well-posedness, long-time behavior, and so on) on more general liquid crystal systems with stretching terms (see also \cite{CRW,WXL13} for the full Ericksen--Leslie system).
 Whether the results obtained in this paper can be extended to those nonautonomous general liquid crystal systems involving stretching effect remains a challenging open problem.
\er

  In order to obtain proper energy inequalities for the system \eqref{1}--\eqref{5}, we recall
that suitable lifting functions were introduced in \cite{C06,C09} to
overcome the technical difficulties related to the time-dependent
boundary datum for $\mathbf{d}$. The first lifting function
$\mathbf{d}_E=\mathbf{d}_E(x,t)$ is of elliptic type (cf.
\cite{C09}):
 \be
 \begin{cases}
 -\Delta \mathbf{d}_E=\mathbf{0},\qquad \text{ in } \Omega\times \mathbb{R}^+,\\
 \mathbf{d}_E=\mathbf{h},\qquad \text{ on } \Gamma\times \mathbb{R}^+.
 \end{cases}\label{LE}
 \ee
In particular, we  define the lifting function $\mathbf{d}_{E0}$ for
the initial datum:
 \be
 \begin{cases}
 -\Delta \mathbf{d}_{E0}=\mathbf{0},\qquad \text{ in } \Omega,\\
 \mathbf{d}_{E0}=\mathbf{d}_0,\qquad  \text{ on } \Gamma.
 \end{cases}\label{iE0}
 \ee
 Set now
 \be
 \wha{\mathbf{d}}=\mathbf{d}-\mathbf{d}_E.\label{dhat}
 \ee
 Then system \eqref{1}--\eqref{5} can be rewritten into the following form:
  \bea
 \mathbf{v}_t+\mathbf{v}\cdot\nabla \mathbf{v}-\nu \Delta \mathbf{v}
 +\nabla \pi&=&-\Delta \wha{\mathbf{d}}\cdot \nabla \mathbf{d}+\mathbf{g}(t),\label{1E}\\
 \nabla \cdot \mathbf{v} &=& 0,\label{2E}\\
 \wha{\mathbf{d}}_t+\mathbf{v}\cdot\nabla \mathbf{d}
 &=&\Delta \wha{\mathbf{d}}-\mathbf{f}(\mathbf{d})-\p_t\mathbf{d}_E(t)\label{3E}
 \eea
 with homogeneous Dirichlet boundary conditions  and initial conditions
 \bea
 && \mathbf{v} =\mathbf{0},\quad \wha{\mathbf{d}} = \mathbf{0},\qquad \text{ on } \Gamma\times
 \mathbb{R}^+,
 \label{4E}\\
 &&
 \mathbf{v}|_{t=0}=\mathbf{v}_0 ,\quad
 \wha{\mathbf{d}}|_{t=0}=\mathbf{d}_0-\mathbf{d}_{E0},\qquad \text{ in }\Omega.\label{5E}
 \eea
 Note that we have used the well-known identity
 $\nabla \cdot(\nabla \mathbf{d}\odot\nabla \mathbf{d})=\frac12\nabla \left(|\nabla \mathbf{d}|^2\right)
 + \Delta \mathbf{d}\cdot \nabla \mathbf{d}$ to absorb the gradient term into pressure (cf. \cite{LL95}).

Let us introduce the lifted energy
 \be
  \wha{\mathcal{E}}(t)=\frac12\|\mathbf{v}(t)\|^2+\frac12
 \|\nabla\wha{\mathbf{d}}(t)\|^2+\int_\Omega F(\mathbf{d}(t)) dx, \quad t\geq 0.\label{E}
 \ee
 Then we can derive the \emph{basic energy inequality} for system \eqref{1}--\eqref{5}.
  \bl
 \label{BEL} Let the assumptions of Proposition \ref{we} be satisfied for all $T>0$.
 Then, any weak solution which is smooth enough satisfies the following inequality for $t\geq 0$
  \be
 \frac{d}{dt}\wha{\mathcal{E}}(t)
 + \frac{\nu}{2} \|\nabla \mathbf{v}\|^2+\frac12\|\Delta \wha{\mathbf{d}}-\mathbf{f}(\mathbf{d})\|^2
 \leq\frac12 \|\p_t \mathbf{d}_E\|^2+C\|\p_t \mathbf{d}_E\|+ C\|\mathbf{g}\|_{\mathbf{V}^*}^2:=r(t),\label{EN1}
 \ee
 where $C$ is a positive constant independent of $\mathbf{v}$ and $\mathbf{d}$.
 \el
 \begin{proof}
  Multiplying \eqref{1E} and \eqref{3E} by $\mathbf{v}$ and
  $-\Delta \wha{\mathbf{d}}+\mathbf{f}(\mathbf{d})$, respectively,
  integrating over $\Omega$ and adding the results together, we get
 \bea
 && \frac{d}{dt}\left(\frac12\|\mathbf{v}\|^2+\frac12 \|\nabla\wha{\mathbf{d}}\|^2+\int_\Omega F(\mathbf{d}) dx\right)
 +\nu \|\nabla \mathbf{v}\|^2+\|\Delta \wha{\mathbf{d}}-\mathbf{f}(\mathbf{d})\|^2\non\\
 &=& (\p_t\mathbf{d}_E, \Delta \wha{\mathbf{d}})+(\mathbf{g}, \mathbf{v}). \label{EN2}
 \eea
In above, we have used the facts $(\mathbf{v}\cdot \nabla \mathbf{v}, \mathbf{v})=(\nabla P,
\mathbf{v})=(\mathbf{v}\cdot\nabla \mathbf{d},
\mathbf{f}(\mathbf{d}))=0$ due to the impressibility condition $\nabla \cdot\mathbf{v}=0$.
 By the Poincar\'e inequality $\|\mathbf{v}\|\leq C_P\|\nabla \mathbf{v}\|$
 and \eqref{max}, the right-hand side of  \eqref{EN2} can be estimated as follows
 \bea
 &&|(\p_t\mathbf{d}_E, \Delta \wha{\mathbf{d}})+(\mathbf{g}, \mathbf{v})|\non\\
 &\leq & |(\p_t\mathbf{d}_E, \Delta \wha{\mathbf{d}}-\mathbf{f}(\mathbf{d}))|+ |(\p_t\mathbf{d}_E, \mathbf{f}(\mathbf{d}))|+|(\mathbf{g}, \mathbf{v})|\non\\
 &\leq&\|\Delta \wha{\mathbf{d}}-\mathbf{f}(\mathbf{d})\|
 \|\p_t\mathbf{d}_E\|+ \|\mathbf{f}(\mathbf{d})\|\|\p_t\mathbf{d}_E\|+ \|\mathbf{v}\|_{\mathbf{V}}\|\mathbf{g}\|_{\mathbf{V}^*}\non\\
 &\leq& \frac{\nu}{2}\|\nabla\mathbf{v}\|^2+ \frac12\|\Delta \wha{\mathbf{d}}-\mathbf{f}(\mathbf{d})\|^2+ \frac12\|\p_t\mathbf{d}_E\|^2+
 C\|\p_t\mathbf{d}_E\|+C\|\mathbf{g}\|_{\mathbf{V}^*}^2.\non
 \eea
 The proof is complete.
 \end{proof}
 \br
 We fix the calculations in \cite[Lemma 2]{C09} where the term
 $(\partial_t\mathbf{d}_E, \Delta \wha{\mathbf{d}})$ is missing.
 Though it does not affect the proof of existence result, it
 does have influence on the long-time behavior of global solutions (especially on the convergence rate).
 \er

Let us now introduce the following  (Banach) spaces of \textit{translation
bounded} functions
\begin{align*}
L^q_{tb}(0,+\infty;X):=&\Big\{
\mathbf{h}\in L^q_{loc} ([0, +\infty);X): \\
&\ \ \Vert \mathbf{h}\Vert_{L^q_{tb}(0,+\infty;X)}^q:=\ \sup_{t\geq
0}\int_{t}^{t+1}\|\mathbf{h}(\tau)\|^q_{X} d\tau<+\infty\Big\}
\end{align*}
where $X$ is a (real) Banach space and $q\in [1,+\infty)$ is given.

From the basic energy inequality \eqref{EN1}, through a suitable Galerkin approximation scheme,
one can derive  uniform-in-time estimates for
any weak solution (the proof is a minor modification of \cite[Lemma 1.2, Remark
1.1]{B}).

 \bl  \label{lowe} Let the assumptions of Proposition \ref{we} hold for all $T>0$.
 In addition, suppose that
 \begin{align}
 \label{hyp6bisa}
&\mathbf{g}\in L^2(0, +\infty; \mathbf{V}^*),\\
\label{hyp7bisa} &\mathbf{h}\in L^2_{tb}(0,+\infty;\mathbf{H}^{\frac32}(\Gamma)) ,\\
\label{hyp8bisa} &\mathbf{h}_t\in L^2(0,+\infty;
\mathbf{H}^{-\frac12}(\Gamma))\cap L^1(0,+\infty;
\mathbf{H}^{-\frac12}(\Gamma)).
\end{align}
Then a weak solution $(\mathbf{v}, \mathbf{d})$ to problem
 \eqref{1}--\eqref{5} given by Proposition \ref{we} is a global solution on $[0,+\infty)$  and
 fulfills the following uniform  bounds
 \bea
 &&
 \|\mathbf{v}(t)\|\leq C, \quad
 \|\mathbf{d}(t)\|_{\mathbf{H}^1}\leq C, \quad \forall\ t\geq 0, \non\\
 && \int_0^t(\nu\|\nabla \mathbf{v}(\tau)\|^2+ \|(\Delta \mathbf{d}-\mathbf{f}(\mathbf{d}))(\tau)\|^2)
 d\tau \leq C,\quad \forall\ t\geq 0.\non
 \eea
Here $C$ is a positive constant depending on
$\|\mathbf{v}_0\|$, $\|\mathbf{d}_0\|_{\mathbf{H}^1}$,
  $\|\mathbf{g}\|_{L^2(0, +\infty; \mathbf{V}^*)}$,
  $\Vert \mathbf{h}\Vert_{L^2_{tb}(0,+\infty;\mathbf{H}^{\frac32}(\Gamma))}$,
  $\|\mathbf{h}_t\|_{ L^2(0,+\infty;
  \mathbf{H}^{-\frac12}(\Gamma))}$ and
  $\|\mathbf{h}_t\|_{ L^1(0,+\infty; \mathbf{H}^{-\frac12}(\Gamma))}$.
 \el
  Next, we introduce the lifting function $\mathbf{d}_P=\mathbf{d}_P(x,t)$ of parabolic type, which satisfies
 \be
 \begin{cases}
 \p_t\mathbf{d}_P-\Delta \mathbf{d}_P=\mathbf{0},\qquad \text{ in } \Omega\times \mathbb{R}^+,\\
 \mathbf{d}_P=\mathbf{h},\qquad   \text{ on } \Gamma\times \mathbb{R}^+,\\
 \mathbf{d}_P(0)= \mathbf{d}_{E0}, \qquad \text{ in } \Omega.
 \end{cases}\label{LP}
 \ee
  The motivation of introducing the parabolic lifting function $\mathbf{d}_P$ is
 that we now have, by definition,
 $\Delta(\mathbf{d}-\mathbf{d}_P)-\mathbf{f}(\mathbf{d})|_\Gamma=\mathbf{0}$.
 This fact is crucial when we use integration by parts to
 derive some higher-order differential inequalities of
 system \eqref{1}--\eqref{5} (cf.  \cite{C09,
 G09}).
 We note that $\mathbf{d}_P$ in \eqref{LP} is different from the one introduced in
 \cite{C09} as they have different initial values. Both choices are valid for the proof of existence result,
 but the current definition of $\mathbf{d}_P$ is necessary for the study of long-time behavior.
  Denote
  $$\wt{\mathbf{d}}=\mathbf{d}-\mathbf{d}_P.$$
 System \eqref{1}--\eqref{5} can now be rewritten into the following form:
  \bea
 \mathbf{v}_t+\mathbf{v}\cdot\nabla \mathbf{v}-\nu \Delta \mathbf{v}+\nabla P&=&-\Delta \mathbf{d}\cdot \nabla \mathbf{d}+\mathbf{g}(t),\label{1P}\\
 \nabla \cdot \mathbf{v} &=& 0,\label{2P}\\
 \wt{\mathbf{d}}_t+\mathbf{v}\cdot\nabla \mathbf{d}&=&\Delta \wt{\mathbf{d}}-\mathbf{f}(\mathbf{d})\label{3P}
 \eea
 with homogeneous Dirichlet boundary conditions and initial conditions
 \bea
 &&\mathbf{v} =\mathbf{0},\quad \wt{\mathbf{d}} = \mathbf{0},\qquad \text{ on } \Gamma\times
 \mathbb{R}^+,
 \label{4P}
 \\
 &&
 \mathbf{v}|_{t=0}=\mathbf{v}_0,\quad
 \wt{\mathbf{d}}|_{t=0}=\mathbf{d}_0-\mathbf{d}_{E0},\qquad {\rm in }\  \Omega.\label{5P}
 \eea
 In the sequel, we shall frequently use the following lemma (cf. \cite{C09})
 \bl \label{dpe}
 The following equivalence between norms hold
 \bea
 && \|\mathbf{v}\|_{\mathbf{H}^1} \approx \|\nabla \mathbf{v}\|,
 \quad \|\wt{\mathbf{d}}\|_{\mathbf{H}^1}
 \approx \|\nabla \wt{\mathbf{d}}\|,\quad {\rm in}\ \mathbf{H}_0^1(\Omega),\non\\
 && \|\mathbf{v}\|_{\mathbf{H}^2} \approx \|\Delta \mathbf{v}\|,\quad \|\wt{\mathbf{d}}\|_{\mathbf{H}^2}
 \approx \|\Delta \wt{\mathbf{d}}\|,\quad {\rm  in}\ \mathbf{H}_0^1(\Omega)
 \cap \mathbf{H}^2(\Omega),\non\\
 && \|\wt{\mathbf{d}}\|_{\mathbf{H}^3}\approx
 \|\nabla(\Delta\wt{\mathbf{d}})\|+\|\Delta \wt{\mathbf{d}}\|,
 \quad  {\rm in}\ \mathbf{H}_0^1(\Omega)\cap \mathbf{H}^3(\Omega).\non
 \eea
  If $\mathbf{d}$ and $\mathbf{d}_P$ are functions that are smooth enough and $|\mathbf{d}|_{\mathbb{R}^n}\leq 1$, $|\mathbf{d}_P|_{\mathbb{R}^n}\leq 1$, then we have
 \bea
 \|\Delta \mathbf{d}\|&\leq& \|\Delta \mathbf{d}_P\|+\|\Delta \wt{\mathbf{d}}-\mathbf{f}(\mathbf{d})\|+C,\non\\
 \|\nabla \Delta \mathbf{d}\|&\leq& \|\nabla \Delta \mathbf{d}_P\|
 +\|\nabla (\Delta \wt{\mathbf{d}}-\mathbf{f}(\mathbf{d}))\|+C\|\nabla
 \mathbf{d}\|,\non
 \eea
 where $C$ is a positive constant independent of $\mathbf{d}$ and $\mathbf{d}_P$.
 \el
 \noindent Let us introduce the quantity
 $$\mathcal{A}_P(t)=\|\nabla \mathbf{v}(t)\|^2 + \|\Delta \wt{\mathbf{d}}(t)-\mathbf{f}(\mathbf{d}(t))\|^2,
 \quad t\geq 0.$$
 \bl \label{H2D}
 Let $n=2$ and let the assumptions of Lemma \ref{lowe} hold.
 If the weak solution $(\mathbf{v}, \mathbf{d})$ is smooth enough then it satisfies the following inequality
  \be \frac{d}{dt}\mathcal{A}_P(t)\leq C(\mathcal{A}_P^2(t)+\mathcal{A}_P(t) + R_1(t)), \label{2dI}
  \ee
  where
  \be
  R_1(t)=\|\p_t \mathbf{d}_P(t)\|^4+\|\p_t \mathbf{d}_P(t)\|^2
  +\|\nabla \Delta \mathbf{d}_P(t)\|^2+\|\mathbf{g}(t)\|^2.\label{2dr}
  \ee
  Here $C$ is a positive constant depending on $\nu$, $\|\mathbf{v}_0\|$,
  $\|\mathbf{d}_0\|_{\mathbf{H}^1}$,
  $\|\mathbf{g}\|_{L^2(0, +\infty; \mathbf{V}^*)}$,
  $\|\mathbf{h}_t\|_{ L^2(0,+\infty;
  \mathbf{H}^{-\frac12}(\Gamma))}$,
  $\|\mathbf{h}_t\|_{ L^1(0,+\infty; \mathbf{H}^{-\frac12}(\Gamma))}$ and
  $\|\mathbf{h}\|_{L^2_{tb}(0,+\infty;\mathbf{H}^\frac32(\Gamma))}$.
 \el
 \begin{proof}
Taking the time derivative of $\mathcal{A}_P(t)$, we obtain by a
direct calculation that
 \bea
 && \frac12\frac{d}{dt}\mathcal{A}_P(t)+ (\nu \|S \mathbf{v}\|^2+ \|\nabla(\Delta
 \wt{\mathbf{d}}-\mathbf{f}(\mathbf{d})\|^2)\non\\
 &=&  -(S \mathbf{v}, \mathbf{v}\cdot\nabla \mathbf{v})+
 (S\mathbf{v}, \mathbf{g})-(S\mathbf{v}, \Delta \mathbf{d}\cdot \nabla \mathbf{d}) -(\nabla (\mathbf{v}\cdot\nabla \mathbf{d} ), \nabla (\Delta
 \wt{\mathbf{d}}-\mathbf{f}(\mathbf{d})))\non\\
 && -(\mathbf{f}'(\mathbf{d})\mathbf{d}_t, \Delta \wt{\mathbf{d}}-\mathbf{f}(\mathbf{d}))\non\\
 &:=&\sum_{j=1}^{5}I_j.
  \label{2dI1}
 \eea
 To get this identity we have used the fact that
 $\Delta\wt{\mathbf{d}}-\mathbf{f}(\mathbf{d})|_\Gamma=\mathbf{0}$ as well as
 $(S \mathbf{v}, \mathbf{v}_t)=(-\Delta \mathbf{v}, \mathbf{v}_t)$. It is not difficult to see that
 \bea
 |I_1|&\leq& \|S \mathbf{v}\|\|\mathbf{v}\|_{\mathbf{L}^4}
 \|\nabla  \mathbf{v}\|_{\mathbf{L}^4}\non\\
 &\leq&   C\|S \mathbf{v}\|(\|\nabla \mathbf{v}\|^\frac12\|\mathbf{v}\|^\frac12)(\|\Delta
 \mathbf{v}\|^\frac12\|\nabla \mathbf{v}\|^\frac12)\non\\
 &\leq& C\|S \mathbf{v}\|^\frac32\|\nabla \mathbf{v}\|\leq \varepsilon \|S \mathbf{v}\|^2+C\|\nabla
 \mathbf{v}\|^2,\non
 \eea
 \be
 |I_2| \leq
   \varepsilon \|S \mathbf{v}\|^2+C\|\mathbf{g}\|^2.\non
  \ee
  For $I_3$, we have
 \bea
 |I_3|&=&|(S\mathbf{v}, (\Delta \wt{\mathbf{d}}-\mathbf{f}(\mathbf{d}))\cdot \nabla \mathbf{d})
 + (S\mathbf{v}, \mathbf{f}(\mathbf{d})\cdot \nabla \mathbf{d})
 + (S\mathbf{v}, \p_t\mathbf{d}_P\cdot \nabla \mathbf{d})|\non\\
 &\leq& \|S\mathbf{v}\|\|\nabla \mathbf{d}\|_{\mathbf{L}^4}
 \|\Delta \wt{\mathbf{d}}-\mathbf{f}(\mathbf{d})\|_{\mathbf{L^4}}+\|S\mathbf{v}\|\|\nabla \mathbf{d}\|_{\mathbf{L}^\infty}
 \|\p_t \mathbf{d}_P\|\non\\
 &\leq& \varepsilon \|S\mathbf{v}\|^2+ C\|\nabla \mathbf{d}\|^2_{\mathbf{L}^4}
 \|\Delta \wt{\mathbf{d}}-\mathbf{f}(\mathbf{d})\|^2_{\mathbf{L}^4}
 + C\|\nabla \mathbf{d}\|^2_{\mathbf{L}^\infty}\|\p_t \mathbf{d}_P\|^2.
 \non
  \eea
 On account of Lemma \ref{lowe},
  we infer from the Sobolev embedding theorems that
  \bea
  \|\nabla \mathbf{d}\|^2_{\mathbf{L}^4}&\leq& C\|\Delta \mathbf{d}\|\|\nabla \mathbf{d}\|+C \|\nabla \mathbf{d}\|^2\leq
   C\|\Delta \wt{\mathbf{d}}\|+C\|\p_t \mathbf{d}_P\|+C\non\\
  &\leq& C\|\Delta \wt{\mathbf{d}}-\mathbf{f}(\mathbf{d})\|+C\|\p_t
  \mathbf{d}_P\|+C,\non
  \eea
  \bea
  \|\nabla \mathbf{d}\|^2_{\mathbf{L}^\infty}&\leq& C\|\nabla \Delta \mathbf{d}\|\|\nabla \mathbf{d}\|+C\|\nabla \mathbf{d}\|^2\non\\
  &\leq&
   C\|\nabla \Delta \mathbf{d}_P\|+C(1+\|\nabla (\Delta
   \wt{\mathbf{d}}-\mathbf{f}(\mathbf{d}))\|),\non
  \eea
  \bea
  \|\Delta
  \wt{\mathbf{d}}-\mathbf{f}(\mathbf{d})\|^2_{\mathbf{L}^4}&\leq&
  C\|\nabla (\Delta \wt{\mathbf{d}}-\mathbf{f}(\mathbf{d}))\|\|\Delta
  \wt{\mathbf{d}}-\mathbf{f}(\mathbf{d})\|.\non
 \eea
 Using the above estimates, we obtain the estimates for $I_3$ and
 $I_4$:
 \bea
 |I_3|&\leq&   \varepsilon \|S\mathbf{v}\|^2+
 C\|\nabla (\Delta \wt{\mathbf{d}}-\mathbf{f}(\mathbf{d}))\| \|\Delta
 \wt{\mathbf{d}}-\mathbf{f}(\mathbf{d})\|(\|\Delta \wt{\mathbf{d}}-\mathbf{f}(\mathbf{d})\|+\|\p_t
 \mathbf{d}_P\|+1)\non\\
 && + C\|\p_t \mathbf{d}_P\|^2(\|\nabla \Delta \mathbf{d}_P\|+\|\nabla (\Delta \wt{\mathbf{d}}-\mathbf{f}(\mathbf{d}))\|+1)
  \non\\
 &\leq& \varepsilon \|S\mathbf{v}\|^2+ \varepsilon \|\nabla (\Delta \wt{\mathbf{d}}-\mathbf{f}(\mathbf{d}))\|^2
 +C\|\Delta \wt{\mathbf{d}}-\mathbf{f}(\mathbf{d})\|^4+ C\|\Delta \wt{\mathbf{d}}-\mathbf{f}(\mathbf{d})\|^2\non\\
 && + C\|\p_t \mathbf{d}_P\|^2(\|\p_t \mathbf{d}_P\|^2+\|\nabla \Delta
 \mathbf{d}_P\|+1),\non
 \eea
 \bea
 |I_4|&\leq & \|\nabla (\Delta \wt{\mathbf{d}}-\mathbf{f}(\mathbf{d}))\|(\|\nabla \mathbf{v}\|_{\mathbf{L}^4}
 \|\nabla \mathbf{d}\|_{\mathbf{L}^4}+\|\mathbf{v}\|_{\mathbf{L}^\infty}\|\mathbf{d}\|_{\mathbf{H}^2})\non\\
 &\leq& \varepsilon \|\nabla (\Delta \wt{\mathbf{d}}-\mathbf{f}(\mathbf{d}))\|^2+ C\|\nabla \mathbf{v}\|_{\mathbf{L}^4}^2
 \|\nabla \mathbf{d}\|_{\mathbf{L}^4}^2+ C\|\mathbf{v}\|^2_{\mathbf{L}^\infty}(\|\Delta\mathbf{d}\|^2+1)\non\\
 &\leq& \varepsilon \|\nabla (\Delta \wt{\mathbf{d}}-\mathbf{f}(\mathbf{d}))\|^2  +C\|\Delta \mathbf{v}\|\|\nabla \mathbf{v}\|(\|\Delta \wt{\mathbf{d}}-\mathbf{f}(\mathbf{d})\|+\|\p_t \mathbf{d}_P\|+1)\non\\
 && +C\|\Delta \mathbf{v}\|\|\mathbf{v}\|(\|\Delta \wt{\mathbf{d}}-\mathbf{f}(\mathbf{d})\|^2+\|\p_t \mathbf{d}_P\|^2+1)\non\\
 &\leq& \varepsilon \|S\mathbf{v}\|^2+\varepsilon \|\nabla (\Delta \wt{\mathbf{d}}-\mathbf{f}(\mathbf{d}))\|^2 + C\|\nabla \mathbf{v}\|^4+ C\|\Delta \wt{\mathbf{d}}-\mathbf{f}(\mathbf{d})\|^4 \non\\
 && +C\|\nabla \mathbf{v}\|^2+C\|\Delta \wt{\mathbf{d}}-\mathbf{f}(\mathbf{d})\|^2+\|\p_t
 \mathbf{d}_P\|^4.\non
 \eea
We now observe that
  \be
 I_5=-(\mathbf{f}'(\mathbf{d})\wt{\mathbf{d}}_t, \Delta
 \wt{\mathbf{d}}-\mathbf{f}(\mathbf{d}))-
 (\mathbf{f}'(\mathbf{d})\p_t\mathbf{d}_P, \Delta
 \wt{\mathbf{d}}-\mathbf{f}(\mathbf{d})):=I_{5a}+I_{5b}.
 \label{split}
 \ee
Recalling \eqref{3P}, we have
 \bea
 |I_{5a}|&=&|(\mathbf{f}'(\mathbf{d})(\mathbf{v}\cdot\nabla)\mathbf{d}, \Delta \wt{\mathbf{d}}-\mathbf{f}(\mathbf{d}))-(\mathbf{f}'(\mathbf{d})(\Delta \wt{\mathbf{d}}-\mathbf{f}(\mathbf{d})), \Delta \wt{\mathbf{d}}-\mathbf{f}(\mathbf{d}))|\non\\
 &\leq& \|\mathbf{f}'(\mathbf{d})\|_{\mathbf{L}^\infty}(\|\mathbf{v}\|_{\mathbf{L}^4}\|\nabla \mathbf{d}\|_{\mathbf{L}^4}\|\Delta \wt{\mathbf{d}}-\mathbf{f}(\mathbf{d})\|+\|\Delta \wt{\mathbf{d}}-\mathbf{f}(\mathbf{d})\|^2)\non\\
 &\leq& C\|\nabla \mathbf{v}\|\|\mathbf{v}\|\|\nabla \mathbf{d}\|_{\mathbf{L}^4}^2+C\|\Delta \wt{\mathbf{d}}-\mathbf{f}(\mathbf{d})\|^2\non\\
 &\leq& C\|\nabla \mathbf{v}\|^2+C\|\Delta \wt{\mathbf{d}}-\mathbf{f}(\mathbf{d})\|^2+C\|\p_t
 \mathbf{d}_P\|^2,\non
 \eea
 \bea
 |I_{5b}|&\leq& \|\mathbf{f}'(\mathbf{d})\|_{\mathbf{L}^\infty}\|\p_t \mathbf{d}_P\|
 \|\Delta \wt{\mathbf{d}}-\mathbf{f}(\mathbf{d})\|\leq C\|\Delta \wt{\mathbf{d}}-\mathbf{f}(\mathbf{d})\|^2
 +C\|\p_t \mathbf{d}_P\|^2.\label{2dI5b}
 \eea
Finally, collecting the above estimates and taking $\varepsilon$
sufficiently small, we deduce that
 \bea && \frac{d}{dt}\mathcal{A}_P(t)+ (\nu\|S \mathbf{v}\|^2+ \|\nabla (\Delta
 \wt{\mathbf{d}}-\mathbf{f}(\mathbf{d}))\|^2)\non\\
 &\leq& C (\mathcal{A}_P^2(t)+\mathcal{A}_P(t)) + C\|\p_t \mathbf{d}_P\|^2(\|\p_t \mathbf{d}_P\|^2
 +\|\nabla \Delta \mathbf{d}_P\|+1)+\|\mathbf{g}\|^2,\non
  \eea
  which easily implies the inequality \eqref{2dI}.
\end{proof}

Taking advantage of Lemmas \ref{H2D}, \ref{Ap1} and \ref{Ap2},
one can deduce the following results on the regularity of weak
solutions as well as the existence of strong solutions to system
\eqref{1}--\eqref{5} in 2D.
 \bt \label{exe2d}
 Let $n=2$ and let the assumptions of Proposition \ref{we} hold for all $T>0$.
 In addition, suppose that
  \begin{align}
 \label{hyp6bis}
&\mathbf{g}\in L^2(0, +\infty; \mathbf{H}),\\
\label{hyp7bis} &\mathbf{h}\in L^2_{tb}(0, +\infty; \mathbf{H}^{\frac52}(\Gamma)) ,\\
\label{hyp8bis} &\mathbf{h}_t\in L^2 (0,+\infty;
\mathbf{H}^{\frac12}(\Gamma))\cap L^1(0,+\infty;
\mathbf{H}^{-\frac12}(\Gamma)).
 \end{align}

 (i) System \eqref{1}--\eqref{5} admits a unique global weak solution $(\mathbf{v}, \mathbf{d})$
  satisfying
  \bea
 && \|\mathbf{v}(t)\|_{\mathbf{V}}\leq C(1+t^{-1}), \quad
 \|\mathbf{d}(t)\|_{\mathbf{H}^2}\leq C(1+t^{-1}), \quad \forall t>
 0,\non
 \\
 && \int_\delta ^t(\|\mathbf{v}(\tau)\|_{\mathbf{H}^2}^2+\|\mathbf{d}(\tau)\|_{\mathbf{H}^3}^2)d\tau
 \leq C(1+\delta^{-1})T,\quad t \in [\delta,T],\non
 \eea
 where $C$ is a positive constant depending on $\nu$, $\|\mathbf{v}_0\|$,
 $\|\mathbf{d}_0\|_{\mathbf{H}^1}$,
  $\|\mathbf{g}\|_{L^2(0, +\infty; \mathbf{H})}$,
  $\|\mathbf{h}\|_{L^2_{tb}(0,+\infty;\mathbf{H}^\frac52(\Gamma))}$, $\|\mathbf{h}_t\|_{L^2(0,+\infty;
 \mathbf{H}^\frac12(\Gamma))}$,
 $\|\mathbf{h}_t\|_{ L^1(0,+\infty; \mathbf{H}^{-\frac12}(\Gamma))}$.

(ii) If $(\mathbf{v}_0, \mathbf{d}_0)\in \mathbf{V}\times
 \mathbf{H}^2(\Omega)$, then problem \eqref{1}--\eqref{5}
 admits a unique global strong solution $(\mathbf{v}, \mathbf{d})$  satisfying
 \bea
 &&\|\mathbf{v}(t)\|_{\mathbf{V}}\leq C, \quad
 \|\mathbf{d}(t)\|_{\mathbf{H}^2}\leq C, \quad \forall t\geq 0, \label{globstrong1}
 \\
 &&\int_0^t(\|\mathbf{v}(\tau)\|_{\mathbf{H}^2}^2+\|\mathbf{d}(\tau)\|_{\mathbf{H}^3}^2)d\tau
 \leq CT,\quad t \in [0,T], \label{globstrong2}
 \eea
 where $C$ is a positive constant depending on $\nu$, $\|\mathbf{v}_0\|_\mathbf{V}$,
 $\|\mathbf{d}_0\|_{\mathbf{H}^2}$,
  $\|\mathbf{g}\|_{L^2(0, +\infty; \mathbf{H})}$,
  $\|\mathbf{h}\|_{L^2_{tb}(0,+\infty;\mathbf{H}^\frac52(\Gamma))}$, $\|\mathbf{h}_t\|_{L^2(0,+\infty;
 \mathbf{H}^\frac12(\Gamma))}$,
 $\|\mathbf{h}_t\|_{ L^1(0,+\infty; \mathbf{H}^{-\frac12}(\Gamma))}$.
 \et

\begin{remark}
Lemma \ref{lowe} and Theorem \ref{exe2d} still hold when $\mathbf{g}$ and $\mathbf{h}_t$
are translation bounded with respect to time (see \cite{B}).
\end{remark}

Next, we consider the 3D case. Instead of Lemma \ref{H2D}, we have
the following higher-order energy inequality

\bl \label{H3D}  Let $n=3$ and let the assumptions of Lemma \ref{lowe} hold.
If a weak solution $(\mathbf{v}, \mathbf{d})$ is smooth enough then it satisfies the following inequality
\bea &&
\frac{d}{dt}\wt{\mathcal{A}}_P(t)+\left(\nu-c_1\wt{\mathcal{A}}_P(t)\right)
\|S \mathbf{v}\|^2
  + \left(1-\frac{c_2}{\nu^\frac12}\wt{\mathcal{A}}_P(t)\right)\|\nabla (\Delta
 \wt{\mathbf{d}}-\mathbf{f}(\mathbf{d}))\|^2\non\\
 &\leq& C(1+\nu^{-2})( \mathcal{A}_P(t) +R_2(t)),\quad t\geq 0,
   \label{3dIA}
 \eea
  where $\wt{\mathcal{A}}_P(t)=\mathcal{A}_P(t)+1$ and
  \be
  R_2(t)=\|\p_t \mathbf{d}_P(t)\|^2+ \|\p_t \mathbf{d}_P(t)\|^6
  +\|\nabla \p_t \mathbf{d}_P(t)\|^2+\|\mathbf{g}(t)\|^2.\label{3dr}
  \ee
Here $c_1,c_2, C$ are positive constants that may depend on
  $\|\mathbf{v}_0\|$, $\|\mathbf{d}_0\|_{\mathbf{H}^1}$ and on
  $\|\mathbf{g}\|_{L^2(0, +\infty; \mathbf{V}^*)}$,
  $\|\mathbf{h}\|_{L^2_{tb}(0,+\infty;\mathbf{H}^\frac32(\Gamma))}$,
  $\|\mathbf{h}_t\|_{ L^2(0,+\infty; \mathbf{H}^{-\frac12}(\Gamma))}$,
   $\|\mathbf{h}_t\|_{ L^1(0,+\infty; \mathbf{H}^{-\frac12}(\Gamma))}$,
  but they are independent of $\nu$.
 \el
 \begin{proof}
 We estimate the right-hand side of \eqref{2dI1} by using
 the 3D version of Sobolev embedding theorems. We have
 \bea
 |I_1|&\leq& \|S \mathbf{v}\|\|\mathbf{v}\|_{\mathbf{L}^6}\|\nabla  \mathbf{v}\|_{\mathbf{L}^3}
 \non\\
 &\leq&  C\|S \mathbf{v}\|\|\nabla \mathbf{v}\|(\|\Delta
 \mathbf{v}\|^\frac12\|\nabla \mathbf{v}\|^\frac12)\non\\
 &\leq& C\|S \mathbf{v}\|^\frac32\|\nabla \mathbf{v}\|^\frac32
 \leq \frac12\|\nabla \mathbf{v}\|^\frac43\|S \mathbf{v}\|^2+C\|\nabla \mathbf{v}\|^2,
 \non
 \eea
 \be |I_2|\leq
   \frac{\nu}{8} \|S
   \mathbf{v}\|^2+\frac{2}{\nu}\|\mathbf{g}\|^2.\non
  \ee
   Recalling that $\|\mathbf{d}\|_{\mathbf{H}^1}\leq C$ (cf. Lemma \ref{lowe}),
   from the Sobolev embedding theorems as well as Agmon's inequality
   in dimension three, we infer
  \be
  \|\nabla \mathbf{d}\|_{\mathbf{L}^3} \leq  C\|\Delta \mathbf{d}\|^\frac12\|\nabla \mathbf{d}\|^\frac12+C \|\nabla \mathbf{d}\|\leq
   C(\|\Delta \wt{\mathbf{d}}-\mathbf{f}(\mathbf{d})\|+\|\p_t \mathbf{d}_P\|)^\frac12+C,
  \non
  \ee
  \be
  \|\nabla \mathbf{d}\|_{\mathbf{L}^6} \leq C\|\Delta \mathbf{d}\|+C \|\nabla \mathbf{d}\|\leq
   C\|\Delta \wt{\mathbf{d}}-\mathbf{f}(\mathbf{d})\|+C\|\p_t \mathbf{d}_P\|+C,
  \non\ee
  \bea
  \|\nabla \mathbf{d}\|_{\mathbf{L}^\infty}
  &\leq& C\|\nabla \mathbf{d}\|_{\mathbf{H}^1}^\frac12\|\nabla \mathbf{d}\|_{\mathbf{H}^2}^\frac12\leq
   C(\|\nabla \Delta \mathbf{d}\|^\frac12\|\Delta \mathbf{d}\|^\frac12+\|\Delta \mathbf{d}\|+1)\non\\
  &\leq& C(\|\nabla (\Delta \wt{\mathbf{d}}-\mathbf{f}(\mathbf{d}))\|^\frac12\|\|\Delta \wt{\mathbf{d}}-\mathbf{f}(\mathbf{d})\|^\frac12+\|\Delta \wt{\mathbf{d}}-\mathbf{f}(\mathbf{d})\|^\frac12\|\nabla \Delta \mathbf{d}_P\|^\frac12\non\\
  &&+\|\nabla (\Delta \wt{\mathbf{d}}-\mathbf{f}(\mathbf{d}))\|^\frac12 \|\Delta \mathbf{d}_P\|^\frac12+\|\nabla \Delta \mathbf{d}_P\|^\frac12\|\Delta \mathbf{d}_P\|^\frac12+ \|\nabla (\Delta \wt{\mathbf{d}}-\mathbf{f}(\mathbf{d}))\|^\frac12\non\\
  && + \|\nabla \Delta \mathbf{d}_P\|^\frac12+\|\Delta \wt{\mathbf{d}}-\mathbf{f}(\mathbf{d})\|+\|\Delta
  \mathbf{d}_P\|+1),\non
  \eea
  \be
  \|\Delta
  \wt{\mathbf{d}}-\mathbf{f}(\mathbf{d})\|_{\mathbf{L}^3} \leq
  C\|\nabla (\Delta \wt{\mathbf{d}}-\mathbf{f}(\mathbf{d}))\|^\frac12
  \|\Delta \wt{\mathbf{d}}-\mathbf{f}(\mathbf{d})\|^\frac12.\non
 \ee
 Thus we have
 \bea
 |I_3| &\leq& \|S\mathbf{v}\|\|\nabla \mathbf{d}\|_{\mathbf{L}^6}
 \|\Delta \wt{\mathbf{d}}-\mathbf{f}(\mathbf{d})\|_{\mathbf{L}^3}
 +\|S\mathbf{v}\|\|\nabla \mathbf{d}\|_{\mathbf{L}^\infty}\|\p_t \mathbf{d}_P\|\non\\
 &\leq& C\|S\mathbf{v}\|(\|\Delta \wt{\mathbf{d}}-\mathbf{f}(\mathbf{d})\|+\|\p_t
 \mathbf{d}_P\|+1)\|\nabla (\Delta \wt{\mathbf{d}}-\mathbf{f}(\mathbf{d}))\|^\frac12
  \|\Delta \wt{\mathbf{d}}-\mathbf{f}(\mathbf{d})\|^\frac12\non\\
  &&+ C\|S\mathbf{v}\|\|\p_t\mathbf{d}_P\|\left(\|\nabla (\Delta \wt{\mathbf{d}}-\mathbf{f}(\mathbf{d}))\|^\frac12\|\Delta \wt{\mathbf{d}}-\mathbf{f}(\mathbf{d})\|^\frac12+\|\Delta \wt{\mathbf{d}}-\mathbf{f}(\mathbf{d})\|^\frac12\|\nabla \p_t\mathbf{d}_P\|^\frac12\right.\non\\
  &&+\|\nabla (\Delta \wt{\mathbf{d}}-\mathbf{f}(\mathbf{d}))\|^\frac12 \|\p_t\mathbf{d}_P\|^\frac12+\|\nabla \p_t\mathbf{d}_P\|^\frac12\|\p_t\mathbf{d}_P\|^\frac12+ \|\nabla (\Delta \wt{\mathbf{d}}-\mathbf{f}(\mathbf{d}))\|^\frac12\non\\
  && \left.+ \|\nabla \p_t\mathbf{d}_P\|^\frac12+\|\Delta \wt{\mathbf{d}}-\mathbf{f}(\mathbf{d})\|+\|\p_t\mathbf{d}_P\|+1\right)\non\\
  &\leq& \left(\frac{\nu}{8}+\frac12 \|\Delta \wt{\mathbf{d}}-\mathbf{f}(\mathbf{d})\|^2\right)\|S\mathbf{v}\|^2
  +\frac18\|\nabla (\Delta \wt{\mathbf{d}}-\mathbf{f}(\mathbf{d}))\|^2\non\\
  && +C (1+\nu^{-2})(
  \|\Delta \wt{\mathbf{d}}
  -\mathbf{f}(\mathbf{d})\|^2
  +\|\p_t \mathbf{d}_P\|^6+\|\p_t\mathbf{d}_P\|^2
  + \|\nabla \Delta\mathbf{d}_P\|^2), \non
  \eea
\bea
 |I_4|&\leq & \|\nabla (\Delta \wt{\mathbf{d}}-\mathbf{f}(\mathbf{d}))\|(\|\nabla \mathbf{v}\|_{\mathbf{L}^3}
 \|\nabla \mathbf{d}\|_{\mathbf{L}^6}+\|\mathbf{v}\|_{\mathbf{L}^\infty}\|\mathbf{d}\|_{\mathbf{H}^2})\non\\
 &\leq& C\|\nabla (\Delta \wt{\mathbf{d}}-\mathbf{f}(\mathbf{d}))\|\|\nabla \mathbf{v}\|^\frac12\|\Delta \mathbf{v}\|^\frac12(\|\Delta \wt{\mathbf{d}}-\mathbf{f}(\mathbf{d})\|+\|\p_t\mathbf{d}_P\|+1)\non\\
 &\leq & \left(\frac{\nu}{8}+\frac12\| \nabla \mathbf{v}\|^2\right)\|S\mathbf{v}\|^2
 +\left(\frac18+\frac{1}{2\nu^\frac12}\|\Delta \wt{\mathbf{d}}-\mathbf{f}(\mathbf{d})\|^2\right)
 \|\nabla (\Delta \wt{\mathbf{d}}-\mathbf{f}(\mathbf{d}))\|^2\non\\
 && +C \|\Delta \wt{\mathbf{d}}-\mathbf{f}(\mathbf{d})\|^2
 +C(1+\nu^{-1})\|\nabla\mathbf{v}\|^2+C\|\p_t\mathbf{d}_P\|^4,\non
  \eea
  \bea
 |I_{5a}| &\leq& \|\mathbf{f}'(\mathbf{d})\|_{\mathbf{L}^\infty}(\|\mathbf{v}\|_{\mathbf{L}^6}\|\nabla \mathbf{d}\|_{\mathbf{L}^3}\|\Delta \wt{\mathbf{d}}-\mathbf{f}(\mathbf{d})\|+\|\Delta \wt{\mathbf{d}}-\mathbf{f}(\mathbf{d})\|^2)\non\\
 &\leq& C\|\nabla \mathbf{v}\|^2\|\nabla \mathbf{d}\|_{\mathbf{L}^3}^2+C\|\Delta \wt{\mathbf{d}}-\mathbf{f}(\mathbf{d})\|^2\non\\
 &\leq& C\|\Delta \mathbf{v}\|\|\mathbf{v}\|(\|\Delta \wt{\mathbf{d}}-\mathbf{f}(\mathbf{d})\|+\|\p_t \mathbf{d}_P\|+1)+C\|\Delta \wt{\mathbf{d}}-\mathbf{f}(\mathbf{d})\|^2\non\\
 &\leq& \frac{\nu}{8}\|S\mathbf{v}\|^2+C(1+\nu^{-1})
 ( \|\nabla \mathbf{v}\|^2+\|\Delta \wt{\mathbf{d}}-\mathbf{f}(\mathbf{d})\|^2+\|\p_t
 \mathbf{d}_P\|^2).\non
 \eea
 We observe that $I_{5b}$ can be estimated as in \eqref{2dI5b}.
 Then, collecting all the estimates of $I_{j}$, we have
 \bea && \frac{d}{dt}\mathcal{A}_P(t)
 + \left(\nu-\|\nabla \mathbf{v}\|^\frac43-\|\nabla \mathbf{v}\|^2\right)
 \|S \mathbf{v}\|^2\non\\
 && + \left(1-\frac{1}{\nu^\frac12}\|\Delta \wt{\mathbf{d}}-\mathbf{f}(\mathbf{d})\|^2\right)\|\nabla (\Delta
 \wt{\mathbf{d}}-\mathbf{f}(\mathbf{d}))\|^2\non\\
 &\leq&  C(1+\nu^{-2})( \mathcal{A}_P(t) + \|\p_t \mathbf{d}_P\|^2+ \|\p_t \mathbf{d}_P\|^6
 +\|\nabla \Delta \mathbf{d}_P\|^2+\|\mathbf{g}\|^2).\non
  \eea
  As a result, there exist constants $c_1,c_2>0$ independent of $\nu$ such that the following inequality holds
  \bea
  && \frac{d}{dt}\mathcal{A}_P(t)+(\nu-c_1\wt{\mathcal{A}}_P(t)) \|S \mathbf{v}\|^2
  + \left(1-\frac{c_2}{\nu^\frac12}\wt{\mathcal{A}}_P(t)\right)\|\nabla (\Delta
 \wt{\mathbf{d}}-\mathbf{f}(\mathbf{d}))\|^2\non\\
 &\leq& C(1+\nu^{-2})(\mathcal{A}_P(t) +\|\p_t \mathbf{d}_P\|^2+ \|\p_t \mathbf{d}_P\|^6+\|\nabla \Delta \mathbf{d}_P\|^2
 +\|\mathbf{g}\|^2),\non
  \eea
  which implies \eqref{3dIA}.
 \end{proof}

On account of Lemma \ref{H3D}, one can deduce that system
\eqref{1}--\eqref{5} admits at least one global strong solution,
provided that the viscosity is large enough (see \cite[Theorem
7]{C09} for the case $\mathbf{g}=\mathbf{0}$, cf. also
\cite{LL95,W10} for the autonomous case). We just report a result
under weaker assumptions than that in \cite{C09} and omit the detailed
proof.

\bt\label{exe3d}
 Let $n=3$ and assume that \eqref{hyp6bis}--\eqref{hyp8bis} and
 \eqref{hyp4} are satisfied.
For any $(\mathbf{v}_0, \mathbf{d}_0)\in \mathbf{V}\times
 \mathbf{H}^2(\Omega)$ satisfying \eqref{hyp5} and
 $|\mathbf{d}_0|_{\mathbb{R}^3}\leq 1$,
 there exists a $\nu_0>0$, depending on $\Vert(\mathbf{v}_0,
\mathbf{d}_0)\Vert_{\mathbf{V}\times
 \mathbf{H}^2}$ and $\Vert\mathbf{g}\Vert_{L^2(0, +\infty; \mathbf{H})}$,
 $\|\mathbf{h}\|_{L^2_{tb}(0,+\infty;\mathbf{H}^\frac52(\Gamma))}$,
$\Vert\mathbf{h}_t\Vert_{L^2(0,+\infty;\mathbf{H}^{\frac12}(\Gamma))}$,
$\|\mathbf{h}_t\|_{ L^1(0,+\infty;
\mathbf{H}^{-\frac12}(\Gamma))}$, such that, for any $\nu\geq
\nu_0$, problem \eqref{1}--\eqref{5} admits a  global strong
solution $(\mathbf{v}, \mathbf{d})$ which satisfies the same
uniform estimates as in the 2D case (cf. \eqref{globstrong1} and
\eqref{globstrong2}).
 \et

\br
 When $n=3$, the weak-strong uniqueness result obtained
 in \cite[Theorem 7]{C09} still
 holds in our case. Thus, the global strong solution $(\mathbf{v},
 \mathbf{d})$ obtained in Theorem \ref{exe3d} is unique.
 \er


\section{Extended \L ojasiewicz--Simon type inequality}
\label{Sec4} \setcounter{equation}{0}

For all $\mathbf{d}\in \mathcal{N}:= \{\phi\in \mathbf{H}^1(\Omega):
\phi|_\Gamma =\mathbf{h}_\infty\}$, where $\mathbf{h}_\infty \in
\mathbf{H}^\frac12(\Gamma)$ is given, we consider the functional
 \be E(\mathbf{d})=\frac12\|\nabla
\mathbf{d}\|^2+\int_\Omega F(\mathbf{d}) dx.\label{fE}
 \ee
It is straightforward to verify that
\bl If $\psi\in \mathbf{H}^1(\Omega)$ is a weak solution to the
elliptic problem \be \left\{\begin{array}{c}
-\Delta\phi+\mathbf{f}(\phi)=\mathbf{0},\\
\phi|_\Gamma=\mathbf{h}_\infty,
\end{array}
\right.\label{sta} \ee then $\psi$ is a critical point of the functional
$E(\mathbf{d})$ in $\mathcal{N}$. Conversely, if $\psi$ is a critical point of the
functional  $E(\mathbf{d})$ in $\mathcal{N}$, then $\psi$ is a weak
solution to problem \eqref{sta}. \el

 \br Due to the elliptic regularity theory, if $\mathbf{h}_\infty$ is more regular, then $\psi$ is more regular. For
instance, if $\mathbf{h}_\infty\in \mathbf{H}^\frac32(\Gamma)$,
then $\psi\in \mathbf{H}^2(\Omega)$.
 \er

 Then we have

\bl
 \label{lsa} Suppose that $\psi$ is a critical point of
$E(\mathbf{d})$ in $\mathcal{N}$.
 Then there exist constants $\beta_1>0$, $\theta\in (0, \frac12)$
 depending on $\psi$ such that, for any
 $\mathbf{w}\in \mathcal{N}$ that satisfies
 $\|\mathbf{w}-\psi\|_{\mathbf{H}^1}<\beta_1$, there holds
\be \label{LSineq}\|-\Delta
\mathbf{w}+\mathbf{f}(\mathbf{w})\|_{\mathbf{H}^{-1}}\geq
|E(\mathbf{w})-E(\psi)|^{1-\theta}.
 \ee
 \el
 \br
 The above lemma can be viewed as an extended version of Simon's result \cite{S83}
 for scalar function under $L^2$-norm.
 We can refer to \cite[Chapter 2, Theorem 5.2]{Hu}, in which the vector case
 subject to homogeneous Dirichlet boundary condition was considered.
 We observe that the result can be easily proved by modifying the
argument in \cite{Hu} using a simple transformation (cf. also
\cite[Remark 2.1]{W10}).
  \er

  The \L ojasiewicz--Simon type inequality \eqref{LSineq} only
  applies to proper perturbations of the critical point of energy $E$ in the set $\mathcal{N}$ and it is not
  enough for our evolutionary problem \eqref{1}--\eqref{5}, whose boundary datum
  is time-dependent (not necessary in $\mathcal{N}$). In order to overcome this difficulty, we prove the following extended
  result that also involves the perturbation of boundary:
\bt \label{lsb} Suppose that $\psi$ is a critical point of
$E(\mathbf{d})$ in $\mathcal{N}$. Then there exists a constant
$\beta\in (0,1)$ depending on $\psi$ such that, for any
$\mathbf{d}\in \mathbf{H}^1(\Omega)$ satisfying
$\|\mathbf{d}-\psi\|_{\mathbf{H}^1}<\beta$, there holds
 \be C\left(\|\mathbf{d}|_\Gamma-\mathbf{h}_\infty\|_{\mathbf{H}^\frac12(\Gamma)}
+\|\mathbf{d}|_\Gamma-\mathbf{h}_\infty\|^{1-\theta}_{\mathbf{H}^\frac12(\Gamma)}\right)
+\|-\Delta \mathbf{d}+\mathbf{f}(\mathbf{d})\|_{\mathbf{H}^{-1}}
\geq |E(\mathbf{d})-E(\psi)|^{1-\theta},\label{ls2}
 \ee
where $\theta\in (0, \frac12)$ is the same constant as in Lemma
\ref{lsa}, while $C$ is a positive constant depending on $\psi$.
 \et
\begin{proof}
For any $\mathbf{d}\in \mathbf{H}^1(\Omega)$, we have that $\Delta
\mathbf{d}\in \mathbf{H}^{-1}(\Omega)$.  Then we consider the
elliptic boundary value problem
 \be
 \left\{\begin{array}{c}
 \Delta\mathbf{w}=\Delta \mathbf{d},\\
 \mathbf{w}|_\Gamma=\mathbf{h}_\infty.
 \end{array}
 \right.
 \label{tran}
 \ee
It easily follows from the elliptic regularity theory (cf. e.g.,
\cite[Proposition 5.1.7]{Tay}) that
 \be
 \|\mathbf{w}- \mathbf{d}\|_{\mathbf{H}^1}\leq
 C\|\mathbf{d}|_\Gamma-\mathbf{h}_\infty\|_{\mathbf{H}^{\frac12}(\Gamma)},\label{we1}
 \ee
which implies
 \bea
 \|\mathbf{w}- \psi\|_{\mathbf{H}^1}&\leq& \|\mathbf{w}- \mathbf{d}\|_{\mathbf{H}^1} +\| \mathbf{d}-\psi\|_{\mathbf{H}^1}
 \non\\
 &\leq& C\|\mathbf{d}|_\Gamma-\mathbf{h}_\infty\|_{\mathbf{H}^\frac12(\Gamma)}+\| \mathbf{d}-\psi\|_{\mathbf{H}^1}\non\\
 &\leq& C\| \mathbf{d}-\psi\|_{\mathbf{H}^1}.\label{we2}
 \eea
Let $\beta_1$ be the constant in Lemma \ref{lsa}. We infer from the
above inequality that if $\beta\in (0,1)$ is chosen sufficiently small,
then we have $\|\mathbf{w}- \psi\|_{\mathbf{H}^1}<\beta_1$.
As a consequence of  Lemma \ref{lsa}, we have
 \be \|-\Delta
\mathbf{w}+\mathbf{f}(\mathbf{w})\|_{\mathbf{H}^{-1}}\geq
|E(\mathbf{w})-E(\psi)|^{1-\theta}.\label{lsw}
 \ee
On the other hand, by the definition of $\mathbf{w}$, we can see
that
 \bea
|E(\mathbf{w})-E(\psi)|^{1-\theta}
 &\leq&  \|-\Delta \mathbf{w}+\mathbf{f}(\mathbf{w})\|_{\mathbf{H}^{-1}}\non\\
 & \leq & \|-\Delta \mathbf{d}+\mathbf{f}(\mathbf{d})\|_{\mathbf{H}^{-1}}
 +C\|\mathbf{f}(\mathbf{d})-\mathbf{f}(\mathbf{w})\|_{\mathbf{L}^\frac65(\Omega)}\non\\
&\leq& \|-\Delta \mathbf{d}+\mathbf{f}(\mathbf{d})\|_{\mathbf{H}^{-1}}+C\|\mathbf{d}-\mathbf{w}\|_{\mathbf{H}^1}\non\\
&\leq& \|-\Delta
\mathbf{d}+\mathbf{f}(\mathbf{d})\|_{\mathbf{H}^{-1}}+C\|\mathbf{d}|_\Gamma-\mathbf{h}_\infty\|_{\mathbf{H}^\frac12(\Gamma)}.
\label{s2}
 \eea
We deduce from $\theta\in (0, \frac12)$ that
 \be
 |E(\mathbf{d})-E(\psi)|^{1-\theta}\leq
|E(\mathbf{w})-E(\psi)|^{1-\theta}
+|E(\mathbf{d})-E(\mathbf{w})|^{1-\theta},\label{s1}
 \ee
and
 \bea && |E(\mathbf{d})-E(\mathbf{w})|^{1-\theta}\non\\
&\leq& \left(\frac12\right)^{1-\theta}\left|\|\nabla\mathbf{d}\|^2- \|\nabla\mathbf{w}\|^2\right|^{1-\theta}+\left|\int_\Omega (F(\mathbf{d})-F(\mathbf{w})) dx\right|^{1-\theta}\non\\
&\leq& C(\|\mathbf{d}\|_{\mathbf{H}^1},
\|\mathbf{w}\|_{\mathbf{H}^1})\|\mathbf{d}-\mathbf{w}\|_{\mathbf{H}^1}^{1-\theta}\leq
C\|\mathbf{d}|_\Gamma-\mathbf{h}_\infty\|^{1-\theta}_{\mathbf{H}^\frac12(\Gamma)},\label{s3}
 \eea
 where in \eqref{s3} we use the facts that $\|\mathbf{d}\|_{\mathbf{H}^1}\leq
 \|\psi\|_{\mathbf{H}^1}+\beta$ and $\|\mathbf{w}\|_{\mathbf{H}^1}\leq
 \|\psi\|_{\mathbf{H}^1}+\beta_1$.
Combining \eqref{s2}--\eqref{s3}, we deduce \eqref{ls2}.
\end{proof}

Since the basic energy inequality \eqref{EN1} (cf. Lemma \ref{BEL})
is only valid for the lifted energy $\wha{\mathcal{E}}$ \eqref{E},
in order to apply the \L ojasiewicz--Simon approach to our problem,
we need to consider the following auxiliary functional corresponding
to energy $E$ (cf. \eqref{fE}):
 \be
 \wha{E}(\mathbf{d})=\frac12\|\nabla \wha{\mathbf{d}}\|^2+\int_\Omega F(\mathbf{d}) dx,
 \quad \forall \ \mathbf{d}\in \mathbf{H}^1(\Omega),
 \label{fEa}
 \ee
 where
 \be
 \wha{\mathbf{d}}=\mathbf{d}-\mathbf{d}_E,\non
 \ee
 and $\mathbf{d}_E$ is the elliptic lifting function satisfying the
 following elliptic problem
 (cf. \eqref{LE})
 \be
 \begin{cases}
 -\Delta \mathbf{d}_E=\mathbf{0},\ \ x\in \Omega,\\
 \mathbf{d}_E=\mathbf{d}|_\Gamma,\quad\  x\in \Gamma.
 \end{cases}\label{aLE}
 \ee
 Then we have
\bc \label{ELS} Suppose that $\psi$ is a critical point of
$E(\mathbf{d})$ in $\mathcal{N}$. Then there exist constants
$\beta\in (0,1)$ and $\theta\in (0,\frac12)$ depending on $\psi$
such that, for any $\mathbf{d}\in \mathbf{H}^1(\Omega)$ satisfying
$\|\mathbf{d}-\psi\|_{\mathbf{H}^1}<\beta$, there holds
 \be C\|\mathbf{d}|_\Gamma-\mathbf{h}_\infty\|^{1-\theta}_{\mathbf{H}^\frac12(\Gamma)}
+\|-\Delta
\wha{\mathbf{d}}+\mathbf{f}(\mathbf{d})\|_{\mathbf{H}^{-1}}\geq
|\wha{E}(\mathbf{d})-\wha{E}(\psi)|^{1-\theta},
 \label{ls4}
 \ee
where $C$ is a positive constant depending on $\psi$ and
$\mathbf{h}_\infty$.
 \ec

\begin{proof}
From the definition of $\wha{E}(\mathbf{d})$, we set, for $\psi\in
\mathcal{N}$,
 \be \wha{E}(\psi)=
 \frac12\|\nabla \wha{\psi}\|^2
 +
 \int_\Omega F(\psi) dx,\label{fEpsi}
 \ee
 where $ \wha{\psi}=\psi-\psi_E$ and $\psi_E$  satisfies
 \be
 \begin{cases}
 -\Delta \psi_E=\mathbf{0},\ \ x\in \Omega,\\
 \psi_E=\mathbf{h}_\infty,\quad\  x\in \Gamma.
 \end{cases}\label{psiLE}
 \ee
 A direct calculation yields that
 \bea
\wha{E}(\mathbf{d})&=& E(\mathbf{d})+\frac12\|\nabla \mathbf{d}_E\|^2-\int_\Omega \nabla \mathbf{d}:\nabla \mathbf{d}_E dx,\non\\
\wha{E}(\psi)&=& E(\psi)+\frac12\|\nabla \psi_E\|^2-\int_\Omega
\nabla \psi:\nabla \psi_E dx,\non
 \eea
 where we used the notation $\mathbf{A}:\mathbf{B}=\sum_{i,j=1}^n A_{ij}B_{ij}$.
 Theorem \ref{lsb} implies
  that there exist constants $\beta\in (0, 1)$ and $\theta\in (0, \frac12)$,
  such that for any $\mathbf{d}\in \mathbf{H}^1(\Omega)$ satisfying $\|\mathbf{d}-\psi\|_{\mathbf{H}^1}<\beta$, \eqref{ls2}
  holds.
 Next, we proceed to estimate the quantity
 $|\wha{E}(\mathbf{d})-\wha{E}(\psi)|^{1-\theta}$
 \bea
 && |\wha{E}(\mathbf{d})-\wha{E}(\psi)|^{1-\theta}
 \non\\
 &\leq& |E(\mathbf{d})-E(\psi)|^{1-\theta}+\left(\frac12\right)^{1-\theta}
 \left|\int_\Omega \nabla (\mathbf{d}_E-\psi_E): \nabla (\mathbf{d}_E+\psi_E)dx\right|^{1-\theta}
 \non\\
 && +\left|\int_\Omega (\nabla \mathbf{d}:\nabla \mathbf{d}_E - \nabla \psi:\nabla \psi_E)dx\right|^{1-\theta}\non\\
 &:=&J_1+J_2+J_3.\label{hatEd}
 \eea
  The estimate for $J_1$ follows from \eqref{ls2}. Since $\|\mathbf{d}-\psi\|_{\mathbf{H}^1}<\beta<1$, then
  $ \|\mathbf{d}\|_{\mathbf{H}^1}\leq \|\psi\|_{\mathbf{H}^1}+1$.  For $J_2$, we infer from the elliptic estimate (cf.
\cite[Proposition 5.1.7]{Tay}) that
   \bea
 J_2 &\leq& C\|\nabla (\mathbf{d}_E-\psi_E)\|^{1-\theta}\| \nabla (\mathbf{d}_E+\psi_E)\|^{1-\theta}\non\\
 &\leq& C\|\mathbf{d}|_\Gamma-\mathbf{h}_\infty\|_{\mathbf{H}^\frac12(\Gamma)}^{1-\theta}
 \left(\|\mathbf{d}|_\Gamma\|_{\mathbf{H}^\frac12(\Gamma)}+\|\mathbf{h}_\infty\|_{\mathbf{H}^\frac12(\Gamma)}\right)^{1-\theta}\non\\
 &\leq& C\|\mathbf{d}|_\Gamma-\mathbf{h}_\infty\|_{\mathbf{H}^\frac12(\Gamma)}^{1-\theta}\left(\|\mathbf{d}\|_{\mathbf{H}^1(\Omega)}
 +\|\mathbf{h}_\infty\|_{\mathbf{H}^\frac12(\Gamma)}\right)^{1-\theta}\non\\
 &\leq& C\|\mathbf{d}|_\Gamma-\mathbf{h}_\infty\|_{\mathbf{H}^\frac12(\Gamma)}^{1-\theta}.
 \eea
 Recalling the function $\mathbf{w}$ introduced in
\eqref{tran}, we estimate  $J_3$ as follows
 \bea
 J_3
 &=&\left|\int_\Omega \left[ \nabla (\mathbf{d}-\mathbf{w}):\nabla \mathbf{d}_E
 +\nabla \mathbf{w}: \nabla (\mathbf{d}_E-\psi_E)+\nabla (\mathbf{w}-\psi): \nabla \psi_E\right] dx\right|^{1-\theta}
\non\\
&\leq& \left|\int_\Omega  \nabla (\mathbf{d}-\mathbf{w}):\nabla \mathbf{d}_E dx\right|^{1-\theta}
 +\left|\int_\Omega \nabla \mathbf{w}: \nabla (\mathbf{d}_E-\psi_E)\right|^{1-\theta}\non\\
 && +\left|\int_\Omega \nabla (\mathbf{w}-\psi): \nabla \psi_E dx\right|^{1-\theta}
 \non\\
 &:=& J_{3a}+J_{3b}+J_{3c}.
 \eea
 Using \eqref{we1} and \eqref{we2} and the fact $\|\mathbf{d}-\psi\|_{\mathbf{H}^1}<\beta$, we observe that
 \bea
 J_{3a}&\leq& \|\nabla (\mathbf{d}-\mathbf{w})\|^{1-\theta}\|\nabla \mathbf{d}_E\|^{1-\theta}\non\\
 &\leq&  C\|\mathbf{d}|_\Gamma-\mathbf{h}_\infty\|_{\mathbf{H}^\frac12(\Gamma)}^{1-\theta}
 \|\mathbf{d}|_\Gamma\|_{\mathbf{H}^\frac12(\Gamma)}^{1-\theta}\non\\
 &\leq& C\|\mathbf{d}|_\Gamma-\mathbf{h}_\infty\|_{\mathbf{H}^\frac12(\Gamma)}^{1-\theta}
 \|\mathbf{d}\|_{\mathbf{H}^1(\Omega)}^{1-\theta}\non\\
 &\leq&
 C\|\mathbf{d}|_\Gamma-\mathbf{h}_\infty\|_{\mathbf{H}^\frac12
 (\Gamma)}^{1-\theta},
 \eea
 \bea
 J_{3b}&\leq& \|\nabla \mathbf{w}\|^{1-\theta} \| \nabla
 (\mathbf{d}_E-\psi_E)\|^{1-\theta}\non\\
 &\leq&
  C\left(\|\psi\|_{\mathbf{H}^1}+C\beta\right)^{1-\theta}
 \|\mathbf{d}|_\Gamma-\mathbf{h}_\infty\|_{\mathbf{H}^\frac12(\Gamma)}^{1-\theta}\non\\
 &\leq& C\|\mathbf{d}|_\Gamma-\mathbf{h}_\infty\|_{\mathbf{H}^
 \frac12(\Gamma)}^{1-\theta}.
 \eea
 For $J_{3c}$, using integration by parts and noticing
 that $\Delta \psi_E=\mathbf{0}$, $(\mathbf{w}-\psi)|_\Gamma=\mathbf{0}$, we obtain
 \be
 \int_\Omega \nabla (\mathbf{w}-\psi) : \nabla \psi_E dx
 = -\int_\Omega (\mathbf{w}-\psi) \cdot \Delta \psi_E dx
   + \int_\Gamma (\mathbf{w}-\psi)|_\Gamma \cdot \partial_\mathbf{n}\psi_E dS=0,\label{lsi3}
 \ee
 where $\mathbf{n}$ is the unit outer  normal to the boundary $\Gamma$. Thus
 \eqref{lsi3} implies that
 \be
 J_{3c}=0. \label{lsi3a}
 \ee
 Finally, since $1-\theta\in(0,1)$, we have
 $\|\mathbf{d}|_\Gamma-\mathbf{h}_\infty\|_{\mathbf{H}^\frac12(\Gamma)}\leq
 C\|\mathbf{d}|_\Gamma-\mathbf{h}_\infty\|_{\mathbf{H}^\frac12(\Gamma)}^{1-\theta}$.
 In summary, we can conclude from \eqref{ls2}, \eqref{hatEd}--\eqref{lsi3a},
 and $\Delta \wha{\mathbf{d}}=\Delta\mathbf{d}$ that \eqref{ls4} holds.
 The proof is complete.
\end{proof}
 \br\label{st}
If $\theta\in (0, \frac12)$ is such that \eqref{ls4} holds, then,
for all $\theta'\in (0, \theta)$ and any $\mathbf{d}\in
\mathbf{H}^1(\Omega)$ satisfying
$\|\mathbf{d}-\psi\|_{\mathbf{H}^1}<\beta$, we still have
 \be C\left(\|\mathbf{d}|_\Gamma-\mathbf{h}_\infty\|^{1-\theta'}_{\mathbf{H}^\frac12(\Gamma)}
+\|-\Delta
\wha{\mathbf{d}}+\mathbf{f}(\mathbf{d})\|_{\mathbf{H}^{-1}}\right)\geq
|\wha{E}(\mathbf{d})-\wha{E}(\psi)|^{1-\theta'},
 \label{ls4a}
 \ee
 where $C$ is a (properly adjusted) positive constant depending on $\psi$
 and $\mathbf{h}_\infty$.
 To see this, we first notice that, since
 $2>\frac{1-\theta'}{1-\theta}>1$, for any $a,b\geq
 0$, it holds
 $(a+b)^\frac{1-\theta'}{1-\theta}\leq
 2(a^\frac{1-\theta'}{1-\theta}+b^\frac{1-\theta'}{1-\theta})$.
 Then it follows from \eqref{ls4}  that
 \bea
 && |\wha{E}(\mathbf{d})-\wha{E}(\psi)|^{1-\theta'}=
 \left(|\wha{E}(\mathbf{d})-\wha{E}(\psi)|^{1-\theta}\right)^\frac{1-\theta'}{1-\theta}\non\\
&\leq&
C^\frac{1-\theta'}{1-\theta}\left(\|\mathbf{d}|_\Gamma-\mathbf{h}_\infty\|^{1-\theta}_{\mathbf{H}^\frac12(\Gamma)}
+\|-\Delta
\wha{\mathbf{d}}+\mathbf{f}(\mathbf{d})\|_{\mathbf{H}^{-1}}\right)^\frac{1-\theta'}{1-\theta}\non\\
&\leq&
C\left(\|\mathbf{d}|_\Gamma-\mathbf{h}_\infty\|^{1-\theta'}_{\mathbf{H}^\frac12(\Gamma)}
+\|-\Delta
\wha{\mathbf{d}}+\mathbf{f}(\mathbf{d})\|_{\mathbf{H}^{-1}}\right).\non
 \eea
 \er

\section{Long-time behavior in 2D}
\label{Sec5} \setcounter{equation}{0}

In this section, we focus on the case $n=2$. In order to study the
long-time behavior of global solutions to problem
\eqref{1}--\eqref{5}, we need some decay conditions on the time-dependent external
force $\mathbf{g}$ and boundary data
$\mathbf{h}$, namely,
\begin{itemize}

\item[(H1)]
    $\int_t^{+\infty}\|\mathbf{h}_t(\tau)\|_{\mathbf{H}^\frac12(\Gamma)}d\tau\leq
    C(1+t)^{-1-\gamma}$;

    \item[(H2)] $\int_t^{+\infty}
        \|\mathbf{h}_t(\tau)\|^2_{\mathbf{H}^\frac12(\Gamma)}
        d\tau\leq C(1+t)^{-1-\gamma}$;

    \item[(H3)] $\int_t^{+\infty}\|\mathbf{g}(\tau)\|^2
        d\tau\leq C(1+t)^{-1-\gamma}$;

 \item[(H4)] $ \|\mathbf{g}(t)\|^2\leq C(1+t)^{-2-\gamma}$;

 \item[(H5)]
     $\|\mathbf{h}_t(t)\|_{\mathbf{L}^2(\Gamma)}\leq
     C(1+t)^{-1-\gamma}$;
\end{itemize}
for all $t\geq 0$. Here $C$ and $\gamma$ are given positive
constants. We also note that (H4) entails (H3).

Since in the 2D case weak solutions become strong for positive times
(cf. Theorem \ref{exe2d}), we can confine ourselves to consider
strong solutions. We recall that, for any given global strong
solution $(\mathbf{v}, \mathbf{d})$, we have the uniform estimate
\eqref{globstrong1}. It follows that the $\omega$-limit set of
the corresponding initial datum $(\mathbf{v}_0,\mathbf{d}_0)$ is
non-empty. Namely, for any unbounded increasing sequence
$\{t_n\}_{n=1}^\infty$, there are functions $\mathbf{v}_\infty\in
\mathbf{V}$ and $\mathbf{d}_\infty\in \mathbf{H}^2(\Omega)$ such
that, up to a subsequence $\{t_j\}_{j=1}^\infty\subset
\{t_n\}_{n=1}^\infty$, we have
 \be \lim_{j\to +\infty}\|\mathbf{v}(t_j)-\mathbf{v}_\infty\|=0,\quad
 \lim_{j\rightarrow +
 \infty}\|\mathbf{d}(t_j)-\mathbf{d}_\infty\|_{\mathbf{H}^1}=0.\label{sconv}
 \ee
 Next, we characterize the structure of the $\omega$-limit set.
In order to do that, we first recall a technical lemma (see \cite[Lemma 6.2.1]{Z04})

\bl \label{SZ}
 Let $T$ be given with $0<T\leq +\infty$. Suppose that $y$ and $h$ are nonnegative continuous functions defined on
 $[0,T]$ and satisfy the following conditions: $
 \frac{dy}{dt}\leq c_1 y^2+ c_2 +h,
 $
 with
 $
 \int_0^T y(t) dt\leq c_3$, $ \int_0^T h(t)dt\leq c_4,
 $
 where $c_i\ (i=1,2,3,4)$ are given nonnegative constants. Then for
 any $\rho\in (0,T)$, the following estimates holds:
 $
 y(t+\rho)\leq \left(\frac{c_3}{\rho}+c_2\rho+c_4\right)e^{c_1c_3}$,
 for all $ t\in[0,T-\rho]$.
 Furthermore, if $T=+\infty$, then
 $\displaystyle\lim_{t\to +\infty} y(t)=0$.
 \el

\bp \label{omega} Let the assumptions of Theorem \ref{exe2d} hold.
Then the
 $\omega$-limit set $\omega(\mathbf{v}_0, \mathbf{d}_0)$ is a subset
 of
 $$\mathcal{S}= \{(\mathbf{0},\mathbf{u}): \mathbf{u}\in \mathcal{N}\cap \mathbf{H}^2(\Omega) \
 {\rm such\ that}\ -\Delta\mathbf{u}+\mathbf{f}(\mathbf{u})=0 \;{ \rm in }\ \Omega\}.$$
Moreover, we have
 \be \lim_{t\to +\infty}\|\mathbf{v}(t)\|_{\mathbf{V}}=0,\label{conv}\ee
 \be \lim_{t\to
 +\infty}\|-\Delta \mathbf{d}(t)+\mathbf{f}(\mathbf{d}(t))\|=0. \label{cond1}\ee
 \ep
\begin{proof}
It follows from Lemma \ref{BEL} that
 $$ \int_0^{+\infty} \|\nabla \mathbf{v}(t)\|^2+\|\Delta
 \wha{\mathbf{d}}(t)-\mathbf{f}(\mathbf{d}(t))\|^2 dt <+\infty,$$
 which together with the definition of $\mathcal{A}_P$ and \eqref{A5} yields
 \be
 \int_0^{+\infty} \mathcal{A}_P(t) dt
 \leq \int_0^{+\infty} (\|\nabla \mathbf{v}(t)\|^2+2\|(\Delta
 \wha{\mathbf{d}}-\mathbf{f}(\mathbf{d}))(t)\|^2+2
 \|\partial_t\mathbf{d}_P(t)\|^2)dt<+\infty.\label{dAP}
 \ee
 Using Lemma \ref{H2D} and Lemma \ref{SZ}, we can see that
 \be \lim_{t\to +\infty } \mathcal{A}_P(t)=0, \non \ee
 which implies $\lim_{t\to +\infty } \|\nabla \mathbf{v}(t)\|=0$. Hence, for any $(\mathbf{v}_\infty,
 \mathbf{d}_\infty)\in \omega(\mathbf{v}_0, \mathbf{d}_0)$, we have $\mathbf{v}_\infty=\mathbf{0}$.
 On the other hand, by definition of $\mathcal{A}_P$, \eqref{dAP} also yields that
 \be
 \lim_{t\to
 +\infty}\|-\Delta \wt{\mathbf{d}}(t)+\mathbf{f}(\mathbf{d}(t))\|=0. \label{cdtt}
 \ee
From Lemma \ref{Ap2}, we have $\displaystyle\lim_{t\to
+\infty}\|\partial_t\mathbf{d}_P(t)\|=0$ (cf. \eqref{dedpt}). As a
 result, it follows from the inequality
 \be
 0\leq \|-\Delta
 \mathbf{d}(t)+\mathbf{f}(\mathbf{d}(t))\|\leq \|-\Delta
 \wt{\mathbf{d}}(t)+\mathbf{f}(\mathbf{d}(t))\|+
 \|\partial_t\mathbf{d}_P(t)\|, \quad \forall\,t\geq 0\label{es1}
  \ee
 that \eqref{cond1} holds. Concerning the limit function $\mathbf{d}_\infty$, we infer from
 \eqref{globstrong1} that $\mathbf{d}_\infty\in \mathbf{H}^2(\Omega)$ and \eqref{sconv} holds.
We now check the boundary condition for $\mathbf{d}_\infty$. Since
$\mathbf{h}_t\in L^1 (0,+\infty; \mathbf{H}^{-\frac12}(\Gamma))$,
$\mathbf{h}(t)$ strongly converges to a certain function
$\mathbf{h}_\infty\in \mathbf{H}^{-\frac12}(\Gamma)$ as time goes to
infinity with a controlled rate, namely,
 \be
 \|\mathbf{h}(t)-\mathbf{h}_\infty\|_{\mathbf{H}^{-\frac12}(\Gamma)}
 \leq\int_t^{+\infty}\|\mathbf{h}_t(\tau)\|_{\mathbf{H}^{-\frac12}(\Gamma)}d\tau\to 0, \quad \text{as} \ t \to +\infty.\label{conh}
 \ee
 On the other hand, we infer from \eqref{hyp7bis} and \eqref{hyp8bis} that $\mathbf{h}\in L^\infty(0,+\infty;
 \mathbf{H}^\frac32(\Gamma))$. Consequently, $\mathbf{h}_\infty\in
 \mathbf{H}^\frac32(\Gamma)$ and $\mathbf{h}$ weakly-star converges
 to $\mathbf{h}_\infty$ in $L^\infty(0,+\infty;
 \mathbf{H}^\frac32(\Gamma))$. By interpolation, we have
 $\displaystyle\lim_{t\to+\infty}\|\mathbf{h}(t)-\mathbf{h}_\infty\|_{\mathbf{L}^2(\Gamma)}=0$.
Thus, from the asymptotic behavior of the
 boundary datum $\mathbf{h}$, we have for any $j\in \mathbb{N}$,
 \bea
 \|\mathbf{d}_\infty|_\Gamma-\mathbf{h}_\infty\|_{\mathbf{L}^2(\Gamma)}
 &\leq&
 \|\mathbf{d}_\infty|_\Gamma-\mathbf{h}(t_j)\|_{\mathbf{L}^2(\Gamma)}+\|\mathbf{h}(t_j)-\mathbf{h}_\infty\|_{\mathbf{L}^2(\Gamma)}\non\\
 &\leq&
 C\|\mathbf{d}_\infty-\mathbf{d}(t_j)\|_{\mathbf{H}^1}+\|\mathbf{h}(t_j)-\mathbf{h}_\infty\|_{\mathbf{L}^2(\Gamma)}.\non
 \eea
Hence, letting $j\to +\infty$ in the above inequality,
 we deduce from \eqref{sconv} and \eqref{conh}  that $\mathbf{d}_\infty|_\Gamma=\mathbf{h}_\infty$.
 For any $\mathbf{z}\in \mathbf{H}_0^1(\Omega)$ and $j\in \mathbb{N}$, we have
 \bea
 && \left|\int_\Omega (-\Delta
 \mathbf{d}_\infty+\mathbf{f}(\mathbf{d}_\infty))\cdot \mathbf{z}dx\right|\non\\
 &\leq& \left|\int_\Omega (-\Delta \mathbf{d}_\infty+\Delta
 \mathbf{d}(t_j))\cdot \mathbf{z} dx\right|+\left|\int_\Omega
 (\mathbf{f}(\mathbf{d}_\infty)-\mathbf{f}(\mathbf{d}(t_j)))\cdot
 \mathbf{z}dx\right|\non\\
 && +\left|\int_\Omega (-\Delta
 \mathbf{d}(t_j)+\mathbf{f}(\mathbf{d}(t_j)))\cdot
 \mathbf{z}dx\right|\non\\
 &\leq&\|\nabla (\mathbf{d}(t_j)-\mathbf{d}_\infty)\|\|\nabla
 \mathbf{z}\|+\big(C\|\mathbf{d}(t_j)-\mathbf{d}_\infty\|_{\mathbf{H}^1}+\|-\Delta
 \mathbf{d}(t_j)+\mathbf{f}(\mathbf{d}(t_j))\|\big)\|\mathbf{z}\|.\non
 \eea
 Passing to the limit as $j\to +\infty$, we get
 \be \int_\Omega (-\Delta
 \mathbf{d}_\infty+\mathbf{f}(\mathbf{d}_\infty))\cdot
 \mathbf{z}dx=0.\non
 \ee
 As a consequence, we see that $\mathbf{d}_\infty \in
 \mathcal{N}\cap\,\mathbf{H}^2(\Omega)$ solves \eqref{sta}. The proof is
 complete.
\end{proof}

We can also prove the convergence of the lifted energy.

 \bp\label{conE}
 Let the assumptions of Theorem \ref{exe2d} hold.
 Then the lifted energy functional
 $\wha{\mathcal{E}}$ defined by \eqref{E}
  is constant on the $\omega$-limit set $\omega(\mathbf{v}_0, \mathbf{d}_0)$. Namely, there exists a constant
  $\wha{\mathcal{E}}_\infty$ such that $\wha{E}(\mathbf{d}_\infty)\equiv\wha{\mathcal{E}}_\infty$, for all $(\mathbf{0},\mathbf{d}_\infty)$ with $\mathbf{d}_\infty \in
 \mathcal{N}\cap\mathbf{H}^2(\Omega)$. Moreover, we have
 \be
\lim_{t\to
+\infty}\wha{\mathcal{E}}(t)=\wha{\mathcal{E}}_\infty.\label{conhE}
 \ee
 \ep
\begin{proof}
From the previous argument, we know that for arbitrary
$(\mathbf{0}, \mathbf{d}^{(1)}_\infty), (\mathbf{0},
\mathbf{d}^{(2)}_\infty)\in \omega(\mathbf{v}_0, \mathbf{d}_0)$
there exist unbounded increasing sequences
$\{t^{(1)}_j\}_{j=1}^\infty$ and $\{t^{(2)}_j\}_{j=1}^\infty$ such
that \eqref{sconv} holds. As a result, we have
 \be\lim_{j\to
 +\infty}\wha{\mathcal{E}}(t^{(1)}_j)=\wha{E}(\mathbf{d}^{(1)}_\infty),\quad
 \lim_{j\to
 +\infty}\wha{\mathcal{E}}(t^{(2)}_j)=\wha{E}(\mathbf{d}^{(2)}_\infty).\non
 \ee
 On the other hand, it follows from the basic energy inequality \eqref{EN1} that
 for any $t'>t''>0$,
 \be |\wha{\mathcal{E}}(t')-\wha{\mathcal{E}}(t'')|\leq
 \int_{t''}^{t'} r(t)dt\to 0, \quad \text{as}\  t', t''\to
 +\infty.\non
  \ee
  Then by
 \be
 |\wha{E}(\mathbf{d}^{(1)}_\infty)-\wha{E}(\mathbf{d}^{(2)}_\infty)|
 \leq
 |\wha{\mathcal{E}}(t_j^{(1)})-\wha{\mathcal{E}}(t_j^{(2)})|+|\wha{E}(\mathbf{d}^{(1)}_\infty)
 -\wha{\mathcal{E}}(t_j^{(1)})|
 +|\wha{\mathcal{E}}(t_j^{(2)})-\wha{E}(\mathbf{d}^{(2)}_\infty)|,\non
 \ee
 letting $j\to +\infty$, we can see that
 $\wha{E}(\mathbf{d}^{(1)}_\infty)=\wha{E}(\mathbf{d}^{(2)}_\infty)$. Namely, $\wha{\mathcal{E}}$
  is a constant (denoted by $\wha{\mathcal{E}}_\infty$) on the $\omega$-limit set $\omega(\mathbf{v}_0,
  \mathbf{d}_0)$.
 Moreover, for any $t>0$ there exist $t_j<t_{j+1}$ such that $t\in
 [t_j,t_{j+1}]$ and
 $|\wha{\mathcal{E}}(t)-\wha{\mathcal{E}}_\infty|\leq
 |\wha{\mathcal{E}}(t)-\wha{\mathcal{E}}(t_j)|+|\wha{\mathcal{E}}(t_j)-\wha{\mathcal{E}}_\infty|$,
 which yields \eqref{conhE}.
\end{proof}

 \subsection{Convergence to equilibrium}

\bt \label{conv1}
 Let the assumptions of Theorem \ref{exe2d} hold.
If, in addition, we assume (H1)--(H3), then
 any strong solution $(\mathbf{v}(t), \mathbf{d}(t))$ convergence
 to an equilibrium $(\mathbf{0}, \mathbf{d}_\infty)$ strongly in
$\mathbf{V}\times \mathbf{H}^2(\Omega)$ as $t$ goes to
$+\infty$. \et

\begin{proof}
On account of \eqref{conv} we only need to prove that
$\mathbf{d}(t)$ converges to $\mathbf{d}_\infty$ as $t \to +\infty$
given by \eqref{sconv}. Below we adapt the idea in \cite{CJ03,GPS1}
to achieve our goal. Indeed, observe that we can find
 an integer $j_0$ such that for all $j\geq j_0$,
 $\|\mathbf{d}(t_j)-\mathbf{d}_\infty\|_{\mathbf{H^1}}<\frac{\beta}{3}$,
 where $\beta\in (0, 1)$ is the constant given in Corollary
 \ref{ELS} (depending on $\mathbf{d}_\infty$). Consequently, we define
 \be
  s(t_j)=\sup\{\tau\geq t_j:
 \|\mathbf{d}(\tau)-\mathbf{d}_\infty\|_{\mathbf{H}^1}<\beta\}.\non
 \ee
Since $\mathbf{d}\in C([0,+\infty); \mathbf{H}^1(\Omega))$, we
 can see that $s(t_j)>t_j$ for any $j \geq j_0$.
 By Lemma \ref{BEL} and Proposition \ref{conE}, we have
 \be
 |\wha{\mathcal{E}}(t)-\wha{E}(\mathbf{d}_\infty)|\geq
 \frac14\min\{\nu,1\}\int_t^{+\infty} \mathcal{D}^2(\tau)d\tau -\int_t^{+\infty} r(\tau)
 d\tau,\non
 \ee
 where
 \be \mathcal{D}(t)=\|\nabla
 \mathbf{v}(t)\|+\|\Delta
 \wha{\mathbf{d}}(t)-\mathbf{f}(\mathbf{d}(t))\|,\non
 \ee
 and $r$ is defined in \eqref{EN1} such that, thanks to (H1)--(H3),
 we have
  \be
 \int_t^{+\infty} r(\tau)
 d\tau\leq C(1+t)^{-1-\gamma}, \quad \forall \ t\geq 0.\non
  \ee
 Let the constant $\theta$ be as in Corollary
 \ref{ELS} (depending on $\mathbf{d}_\infty$).  Using Remark \ref{st}, we can choose $\theta'\in (0, \theta]$ such
 that $\theta'$ also satisfies
 \be
 0<\theta'<\frac{\gamma}{2(1+\gamma)}. \label{theta}
 \ee
 If $\theta$ itself satisfies \eqref{theta}, we just take
 $\theta'=\theta$.
 For any fixed $t_j$ with $j\geq j_0$, we introduce the sets
 \be
 K_j=[t_j,s(t_j)),\quad
 K^{(1)}_j=\left\{t\in K_j:
 \mathcal{D}(t)>(1+t)^{-(1-\theta')(1+\gamma)}\right\},\quad
 K_j^{(2)}=K_j\setminus K_j^{(1)}.\non
 \ee
Consider the following functional on $K_j$
 \be
 \Phi(t)=\wha{\mathcal{E}}(t)-\wha{E}(\mathbf{d}_\infty)+2\int_t^{s(t_j)}
 r(\tau)d\tau, \quad \forall \ t\in K_j.\non
 \ee
 It easily follows that
 \be
 \lim_{j\to+\infty} \Phi(t_j)=0.\label{decayp}
 \ee
 Next, we have
 \bea
 \frac{d}{dt}(|\Phi(t)|^{\theta'}{\rm sgn} \Phi(t))
 &=&\theta'|\Phi(t)|^{\theta'-1}\frac{d}{dt}\Phi(t)\non\\
 &\leq&
 -\frac{\theta'}{4}\min\{\nu, 1\}|\Phi(t)|^{\theta'-1}\mathcal{D}^2(t)\non\\
 &\leq&  0,\label{di}
 \eea
 which implies that the functional $|\Phi(t)|^{\theta'}{\rm sgn}
 \Phi(t)$ is decreasing on $K_j$.
  Keeping in mind that $\theta'\leq \theta$ and $2(1-\theta')>1$, we can apply Corollary \ref{ELS} (cf. also Remark \ref{st})
   to obtain that
 \bea
 |\Phi(t)|^{1-\theta'}
 &\leq& |\wha{\mathcal{E}}(t)-\wha{E}(\mathbf{d}_\infty)|^{1-\theta'}+C\left(\int_t^{+\infty}
 r(\tau)d\tau\right)^{1-\theta'}\non\\
 &\leq& \left(\frac12\right)^{2(1-\theta')}\|\mathbf{v}\|^{2(1-\theta')}
 + C\|\mathbf{h}(t)-\mathbf{h}_\infty\|^{1-\theta'}_{\mathbf{H}^\frac12(\Gamma)}\non\\
 &&
+C\|-\Delta
\wha{\mathbf{d}}+\mathbf{f}(\mathbf{d})\|_{\mathbf{H}^{-1}}+C\left(\int_t^{+\infty}
 r(\tau)d\tau\right)^{1-\theta'}\non\\
 &\leq& C\|\nabla \mathbf{v}\|+C\|-\Delta
\wha{\mathbf{d}}+\mathbf{f}(\mathbf{d})\|+
  C\left(\int_t^{+\infty}\|\mathbf{h}_t(\tau)\|_{\mathbf{H}^\frac12(\Gamma)}d\tau\right)^{1-\theta'}\non\\
  && +C\left(\int_t^{+\infty}
 r(\tau)d\tau\right)^{1-\theta'}\non\\
 &\leq& C\|\nabla \mathbf{v}\| +C\|-\Delta
\wha{\mathbf{d}}+\mathbf{f}(\mathbf{d})\|+ C
(1+t)^{-(1-\theta')(1+\gamma)}.\label{dePhi}
 \eea
 Thus, on $K_j^{(1)}$, we have
 \be
 |\Phi(t)|^{1-\theta'}\leq C\mathcal{D}(t),\non
 \ee
 which together with \eqref{di} yields that on
 $K_j^{(1)}$,
 \be
 -\frac{d}{dt}(|\Phi(t)|^{\theta'}{\rm sgn} \Phi(t))
 \geq  C\mathcal{D}(t).\label{di1}
 \ee
As a consequence, we have
 \bea
  \int_{K_j^{(1)}}\mathcal{D}(t)dt &\leq& -C\int_{K_j} \frac{d}{dt}(|\Phi(t)|^{\theta'}{\rm
 sgn}
 \Phi(t))dt\non\\
 &\leq&
 C(|\Phi(t_j)|^{\theta'}+|\Phi(s(t_j))|^{\theta'})<+\infty,\label{D1}
 \eea
 where $\Phi(s(t_j))=0$ if $s(t_j)=+\infty$. On the other hand, on $K_j^{(2)}$, we have
 \be
 \int_{K_j^{(2)}} \mathcal{D}(t) dt\leq C\int_{t_j}^\infty
  (1+t)^{-(1-\theta')(1+\gamma)}dt =
  \frac{C}{-\gamma\theta'-\theta'+\gamma}(1+t_j)^{\gamma\theta'+\theta'-\gamma}.\label{D2}
 \ee
 Here, we notice that $\gamma\theta'+\theta'-\gamma<0$ due to \eqref{theta}. Then \eqref{D1} and \eqref{D2} imply
 that
 $$
 \int_{K_j}
 \mathcal{D}(t) dt =\int_{K_j^{(1)}}\mathcal{D}(t)dt +\int_{K_j^{(2)}} \mathcal{D}(t)
 dt<+\infty,
 $$
 for any $j$. On the other hand, it follows from \eqref{globstrong1} and
 \eqref{3E}  that
 \bea
 \|\mathbf{d}_t(t)\| &\leq&  \|\mathbf{v}\cdot\nabla \mathbf{d}\|+\|\Delta
 \wha{\mathbf{d}}-\mathbf{f}(\mathbf{d})\| \non\\
 &\leq&  \|\mathbf{v}\|_{\mathbf{L}^4}\|\nabla
 \mathbf{d}\|_{\mathbf{L}^4}+\|\Delta
 \wha{\mathbf{d}}-\mathbf{f}(\mathbf{d})\|
 \non\\
 &\leq&
 C\mathcal{D}(t).\label{vadt}
 \eea
 As a consequence,
 \be
 \int_{K_j} \|\mathbf{d}_t(t)\|dt \leq
 C(|\Phi(t_j)|^{\theta'}+|\Phi(s(t_j))|^{\theta'})+C(1+t_j)^{\gamma\theta'+\theta'-\gamma}.\label{sdt1}
 \ee
 To complete the proof, we show that
  \bp
 \label{jinf} Let the assumptions of Theorem \ref{exe2d} hold.
 Then there exists an integer $j_1\geq j_0$ such that
 $s(t_{j_1})=+\infty$. Thus
 $$\|\mathbf{d}(t)-\mathbf{d}_\infty\|_{\mathbf{H}^1}<\beta,\quad
 \forall\,t \geq t_{j_1}.
 $$
 \ep
 \begin{proof}
 The conclusion follows from a contradiction argument (cf. \cite{J981}).
 Suppose that for any $j\geq j_0$ we have $s(t_j)<+\infty$. Then, by definition, we have
 \be
 \|\mathbf{d}(s(t_j))-\mathbf{d}_\infty\|_{\mathbf{H}^1}=\beta>0.\label{codi}
 \ee
 Besides, it follows from \eqref{sconv}, \eqref{decayp} and \eqref{sdt1} that
 \bea
 \|\mathbf{d}(s(t_j))-\mathbf{d}_\infty\|&\leq&
 \|\mathbf{d}(s(t_j))-\mathbf{d}(t_j)\|+\|\mathbf{d}(t_j)-\mathbf{d}_\infty\|\non\\
 &\leq&
 \int_{t_j}^{s(t_j)} \|\mathbf{d}_t(t)\|dt+
 \|\mathbf{d}(t_j)-\mathbf{d}_\infty\|\to 0, \quad \text{as}\ j\to
 +\infty.\non
 \eea
 Using uniform estimate \eqref{globstrong1} and interpolation inequality, we obtain
  \be
  \|\mathbf{d}(s(t_j))-\mathbf{d}_\infty\|^2_{\mathbf{H}^1}
  \leq
  \|\mathbf{d}(s(t_j))-\mathbf{d}_\infty\|_{\mathbf{H}^2}\|\mathbf{d}(s(t_j))-\mathbf{d}_\infty\|\to
  0,\quad \text{as}\ j\to
 +\infty,\non
 \ee
 which leads a contradiction with \eqref{codi}. The proof is
 complete.
 \end{proof}
 Due to Proposition \ref{jinf}, we have $s(t_{j_1})=+\infty$ for some $j_1\geq j_0$.
 Arguing as above, we can prove
 \be
 \int_{t_{j_1}}^{+\infty}\|\mathbf{d}_t(t)\|dt<+\infty.\non
 \ee
 Thus $\mathbf{d}(t)$ converges in $\mathbf{L}^2$ and recalling
 \eqref{sconv}, by compactness we conclude that
 \be
 \label{convH1}
 \lim_{t\to+\infty}
 \|\mathbf{d}(t)-\mathbf{d}_\infty\|_{\mathbf{H}^1}=0.
 \ee
Finally, observe that
 \bea
 \|\Delta \mathbf{d}(t)-\Delta \mathbf{d}_\infty\|
 &=& \|-\Delta
 \mathbf{d}(t)+\mathbf{f}(\mathbf{d}(t))\|+\|\mathbf{f}(\mathbf{d}(t))
 -\mathbf{f}(\mathbf{d}_\infty)\|\non\\
 &\leq& \|-\Delta
 \mathbf{d}(t)+\mathbf{f}(\mathbf{d}(t))\|
 +C\|\mathbf{d}(t)-\mathbf{d}_\infty\|.\label{eD}
 \eea
 Then \eqref{cond1} and \eqref{convH1} entail that
\be
 \lim_{t\to+\infty}\|\mathbf{d}(t)-\mathbf{d}_\infty\|_{\mathbf{H}^2}
 =0\non
 \ee
 and this finishes the proof of Theorem \ref{conv1}.
 \end{proof}
\subsection{Convergence rate}
  \bt\label{convrate1}
 Let the assumptions of Theorem \ref{exe2d} hold.
If, in addition, we assume (H1)--(H2) and (H4)--(H5), then we have
\be
  \|\mathbf{v}(t)\| +\|\mathbf{d}(t)-
  \mathbf{d}_\infty\|_{\mathbf{H}^1}
  \leq   C(1+t)^{- \frac{\theta'}{1-2\theta'}},\quad
  t\geq 0.\non
  \ee
Moreover, if (H2) and (H5) are replaced by, respectively,
\begin{itemize}
\item[(H6)]
    $\|\mathbf{h}_t(t)\|_{\mathbf{H}^{\frac12}(\Gamma)}\leq
    C(1+t)^{-1-\gamma}$;

\item[(H7)] $\|\mathbf{h}(t) - \mathbf{h}_\infty
    \|_{\mathbf{H}^{\frac32}(\Gamma)}\leq
    C(1+t)^{-1-\gamma}$;
\end{itemize}
the following higher-order estimate holds \be
  \|\mathbf{v}(t)\|_{\mathbf{V}}+\|\mathbf{d}(t)-
  \mathbf{d}_\infty\|_{\mathbf{H}^2}
  \leq   C(1+t)^{- \frac{\theta'}{1-2\theta'}},\quad
  t\geq 0.\non
  \ee

 \et

\begin{proof}
The proof consists of several steps.

 \textbf{Step 1}. $\mathbf{L}^2$-estimate of
 $\mathbf{d}-\mathbf{d}_\infty$. This follows from an argument devised in
 \cite{GPS1}. For the readers' convenience, we sketch the
 proof here.  From the previous argument, we only have to work on the time interval $[t_{j_1}, +\infty)$.
 Denote
  \be
 \Phi(t)=\wha{\mathcal{E}}(t)-\wha{E}(\mathbf{d}_\infty)+2\int_t^{+\infty}
 r(\tau)d\tau.\non
 \ee
 Since
 \be
 \frac{d}{dt}\Phi(t) \leq -\frac{\theta'}{4}\min\{\nu,
 1\}\mathcal{D}^2(t)-r(t)\leq 0,\non
 \ee
 and $\displaystyle\lim_{t\to+\infty}\Phi(t)=0$, we know that $\Phi(t)$ is
 decreasing and $\Phi(t)\geq 0$ for $t\geq t_{j_1}$.

First, if the boundary datum $\mathbf{h}$ and the external force
$\mathbf{g}$ become time-independent in finite time, i.e., there
exists time $T_0$ such that for $t\geq T_0$,
$\mathbf{h}=\mathbf{h}_\infty$ and $\mathbf{g}=\mathbf{0}$. Then the
problem reduces to the autonomous system considered in \cite{W10}.
Thus, below we just assume that either $\mathbf{h}$ or $\mathbf{g}$
does not become time-independent in finite time (namely, the system
will always be non-autonomous). In this case, if there exists
$t^*\geq t_{j_1}$ such that $\Phi(t^*)=0$, then
$\mathcal{D}(t)=r(t)=0$ for all $t\geq t^*$ and this is a
contradiction since $r(t)$ cannot identically vanish from any finite
time on. Therefore, we can suppose
 \be \Phi(t)>0, \quad \forall\ t\geq t_{j_1}.\non\ee

 If the open set $K_{j_1}^{(1)}$ is bounded, then there
 exists $t^*\geq t_{j_1}$ such that $[t^*,+\infty)\subset K_{j_1}^{(2)}$.
 As a result, $\mathcal{D}(t)\leq
 (1+t)^{-(1-\theta')(1+\gamma)}$ and by \eqref{vadt},
  we have
 \be
 \|\mathbf{d}(t)-\mathbf{d}_\infty\|\leq \int_t^{+\infty}
 \|\mathbf{d}_t(\tau)\|d\tau\leq
 \frac{C}{-\gamma\theta'-\theta'+\gamma}(1+t)^{\gamma\theta'+\theta'-\gamma}, \quad \forall\ t\geq
 t^*.\non
 \ee

 Next, we treat the case when the open set $K_{j_1}^{(1)}$ is unbounded.
 There exists a countable family of disjoint open sets $(a_n, b_n)$ such that
 $K_{j_1}^{(1)}=\cup_{n=1}^\infty(a_n,b_n)$. On $K_{j_1}^{(1)}$, recalling \eqref{dePhi}, we
 can see that on any $(a_n,b_n)\subset K_{j_1}^{(1)}$, it holds
 \be
  \frac{d}{dt}\Phi(t)+C\Phi^{2(1-\theta')}(t)\leq 0.\non
 \ee
 As a result, for any $t\in (a_n,b_n)$,
  \be \Phi(t)\leq
 \left[\Phi(a_n)^{2\theta'-1}+C(1-2\theta')(t-a_n)\right]^{-\frac{1}{1-2 \theta'}},\label{j1aaa}
 \ee
 where by the definition of $K_{j_1}^{(1)}$ and \eqref{dePhi} we
 have
 \be
 \Phi(a_n)\leq C\mathcal{D}(a_n)^{\frac{1}{1-\theta'}}+ C
(1+a_n)^{-(1+\gamma)}= C (1+a_n)^{-1-\gamma}.\non
 \ee
  Using the fact $(1+\gamma)(1-2\theta')>1$ (cf. \eqref{theta}), we can take $n^*\in {\mathbb N}$
 sufficiently large such that
 \be \Phi(a_{n^*})^{2\theta'-1}-C(1-2\theta')a_{n^*}\geq
 a_{n^*}^{(1+\gamma)(1-2\theta')}-C(1-2\theta
 ')a_{n^*}\geq 1.\label{j1aaa1}
 \ee
 Therefore, we infer
  \be
  \Phi(t)\leq C(1+t)^{-\frac{1}{1-2\theta'}}, \quad \forall\ t\in (a_{n^*},\infty)\cap K_{j_1}^{(1)}.\non
  \ee
  Similar to \eqref{di1}, we have (since $\Phi(t)>0$)
  \be
  -\frac{d}{dt} \Phi(t)^{\theta'} \geq C\mathcal{D}(t), \quad \forall\ t\in (a_{n^*},\infty)\cap
  K_{j_1}^{(1)}.\non
  \ee
  Due to \eqref{theta}, it follows that $-\gamma\theta'-\theta'+\gamma\geq
 \frac{\theta'}{1-2\theta'}$. Now for any $t>a_{n^*}$, we can conclude that
  \bea
  \|\mathbf{d}(t)-\mathbf{d}_\infty\|&\leq&  \int_t^{+\infty}\|\mathbf{d}_t(\tau)\|d\tau\non\\
  &=&
  \int_{ (t,\infty)\cap
  K_{j_1}^{(1)}}\|\mathbf{d}_t(\tau)\|d\tau+\int_{ (t,\infty)\cap
  K_{j_1}^{(2)}}\|\mathbf{d}_t(\tau)\|d\tau\non\\
  &\leq&  C \int_{ (t,\infty)\cap  K_{j_1}^{(1)}}
  \mathcal{D}(\tau)d\tau+ C\int_t^{+\infty} (1+\tau)^{-(1-\theta')(1+\gamma)}d\tau \non\\
  &\leq& C\Phi(t)^{\theta'}+ C(1+t)^{\gamma\theta'+\theta'-\gamma}\non\\
  & \leq&  C(1+t)^{-\frac{\theta'}{1-2\theta'}}.\non
  \eea
Using  \eqref{globstrong1}, after properly adjusting the constant
$C$, we have
 \be
 \|\mathbf{d}(t)-\mathbf{d}_\infty\|\leq
 C(1+t)^{-\frac{\theta'}{1-2\theta'}}, \quad \forall\ t\geq
 0.\label{rate1}
 \ee

 \textbf{Step 2}. $\mathbf{H}\times \mathbf{H}^1$-estimate.
 It easily from the basic energy inequality \eqref{EN1} that
  \be
  \frac{d}{dt}y(t)
  + \frac{\nu}{2} \|\nabla \mathbf{v}\|^2+\frac12\|\Delta \wha{\mathbf{d}}-\mathbf{f}(\mathbf{d})\|^2
  \leq r(t),\label{ENdyt}
 \ee
 where
 \be
 y(t)= \frac12\|\mathbf{v}(t)\|^2+ \frac12\| \nabla
 (\wha{\mathbf{d}}(t)-\wha{\mathbf{d}}_\infty)
 \|^2+ \int_\Omega[
 F(\mathbf{d})(t)-F(\mathbf{d}_\infty)-\mathbf{f}(\mathbf{d}_\infty)(\mathbf{d}(t)-\mathbf{d}_\infty)]dx.\non
 \ee
 As in \cite{W10}, using  \eqref{globstrong1}, we can show that
 \be
 \left|\int_\Omega[
 F(\mathbf{d})(t)-F(\mathbf{d}_\infty)-\mathbf{f}(\mathbf{d}_\infty)(\mathbf{d}(t)-\mathbf{d}_\infty)]dx\right|\leq
 C \|\mathbf{d}(t)-\mathbf{d}_\infty\|^2.\non
 \ee
 Keeping in mind the definition of lifting functions, we have $\wha{\mathbf{d}}-\wha{\mathbf{d}}_\infty|_{\Gamma}=0$
so that
 \bea
  \| \nabla  (\wha{\mathbf{d}}-\wha{\mathbf{d}}_\infty)\|
 &\leq&
 C\|\Delta (\wha{\mathbf{d}}-\wha{\mathbf{d}}_\infty)\|\non\\
 &\leq&
\|-\Delta
 \wha{\mathbf{d}}+\mathbf{f}(\mathbf{d})\|+C\|\mathbf{f}(\mathbf{d})-\mathbf{f}(\mathbf{d}_\infty)\|\non\\
 &\leq& C\|-\Delta
 \wha{\mathbf{d}}+\mathbf{f}(\mathbf{d})\|+C\|\mathbf{d}(t)-\mathbf{d}_\infty\|,\non
 \eea
 \bea
 \| \nabla  (\mathbf{d}-\mathbf{d}_\infty)\|&\leq&  \| \nabla
(\wha{\mathbf{d}}-\wha{\mathbf{d}}_\infty)\|
+ \| \nabla  (\mathbf{d}_E-\mathbf{d}_\infty)\|\non\\
 &\leq& C\|-\Delta
 \wha{\mathbf{d}}+\mathbf{f}(\mathbf{d})\|
 +C\|\mathbf{d}(t)-\mathbf{d}_\infty\|+C\int_t^{+\infty}
 \|\mathbf{h}_t(\tau)\|_{\mathbf{H}^\frac12(\Gamma)}d\tau.\non
 \eea
 Thus it follows that
 \bea
 y(t)&\geq & \frac12\|\mathbf{v}(t)\|^2+ \frac12\| \nabla
 (\mathbf{d}-\mathbf{d}_\infty)
 \|^2-C\|\mathbf{d}(t)-\mathbf{d}_\infty\|^2\non\\
 &&
 -C\left(\int_t^{+\infty} \|\mathbf{h}_t(\tau)\|_{\mathbf{H}^\frac12(\Gamma)}d\tau\right)^2, \label{ya}\\
 y(t)&\leq& C\|\nabla \mathbf{v}\|^2+  C\|-\Delta
 \wha{\mathbf{d}}+\mathbf{f}(\mathbf{d})\|^2 +C\|\mathbf{d}(t)-\mathbf{d}_\infty\|^2.
 \eea
  Condition \eqref{theta} implies that $ \frac{2\theta'}{1-2\theta'}<\gamma$.
  Then we deduce from \eqref{ENdyt}, \eqref{rate1}, (H4)--(H5) and Lemma \ref{Ap1} that
 \be
 \frac{d}{dt}y(t)+
 \alpha y(t)\leq C(r(t)+\|\mathbf{d}(t)-\mathbf{d}_\infty\|^2)\leq C(1+t)^{- \frac{2\theta'}{1-2\theta'}}, \label{rate2}
  \ee
  where $\alpha>0$ is sufficiently small. The above inequality implies
  that
 \be
 y(t)\leq C(1+t)^{- \frac{2\theta'}{1-2\theta'}},\quad \forall\
 t\geq 0.\label{ratey}
 \ee
 Combining it with \eqref{ya} and recalling (H1), we get
 \bea
 && \|\mathbf{v}(t)\|^2+ \| \mathbf{d}(t)-\mathbf{d}_\infty
 \|_{\mathbf{H}^1}^2\non\\
 &\leq& Cy(t)+C\|\mathbf{d}(t)-\mathbf{d}_\infty\|^2
 +C\left(\int_t^{+\infty}
 \|\mathbf{h}_t(\tau)\|_{\mathbf{H}^\frac12(\Gamma)}d\tau\right)^2\non\\
 &\leq&  C(1+t)^{- \frac{2\theta'}{1-2\theta'}},\quad \forall\
 t\geq 0.\label{rate4}
 \eea

  \textbf{Step 3}. $\mathbf{V}\times \mathbf{H}^2$-estimate.
Taking advantage of the stronger assumptions (H6)--(H7) and
\eqref{rate4}, we now get a higher-order estimate. Observe first
that
 \be
 \|-\Delta
 \wha{\mathbf{d}}+\mathbf{f}(\mathbf{d})\|\leq \|-\Delta
 \wt{\mathbf{d}}+\mathbf{f}(\mathbf{d})\| +\|\Delta
 \mathbf{d}_P\|=\|-\Delta
 \wt{\mathbf{d}}+\mathbf{f}(\mathbf{d})\| +\|\partial_t
 \mathbf{d}_P\|,\non
 \ee
 then we have
 \be
 y(t)\leq C\|\nabla \mathbf{v}\|^2+  C\|-\Delta
 \wt{\mathbf{d}}+\mathbf{f}(\mathbf{d})\|^2 + C\|\partial_t
 \mathbf{d}_P\|^2+C\|\mathbf{d}(t)-\mathbf{d}_\infty\|^2.\non
 \ee
It follows from \eqref{2dI} and \eqref{rate2} that
 \be
 \frac{d}{dt}z(t)+
 \alpha_2 z(t)
 \leq C(1+t)^{- \frac{2\theta'}{1-2\theta'}}+ C(R_1(t)+\|\partial_t
 \mathbf{d}_P(t)\|^2),
  \label{rate3}
  \ee
  where
  \be
  z(t)=y(t)+\alpha_1\mathcal{A}_P(t),
  \ee
 and $\alpha_1$ and $\alpha_2$ are sufficiently small positive constants.
 From the definition
 of $R_1$, \eqref{A8} and the fact $\frac{2\theta'}{1-2\theta'}< 2+2\gamma $ (cf. \eqref{theta}),
 we have
 \be
 R_1(t)+\|\partial_t
 \mathbf{d}_P(t)\|^2\leq C(1+t)^{- \frac{2\theta'}{1-2\theta'}}+
 C\|\nabla \Delta \mathbf{d}_P(t)\|^2.
 \ee
 Hence, from \eqref{rate3} we infer that
  \bea z(t)
 &\leq& z(0)
e^{-\alpha_2 t}+Ce^{-\alpha_2 t}\int_0^t e^{\alpha_2\tau}
\left[C(1+\tau
)^{- \frac{2\theta'}{1-2\theta'}}+\|\nabla \Delta \mathbf{d}_P(\tau)\|^2\right]d\tau\non\\
&\leq& Ce^{-\alpha_2 t}
 +e^{-\alpha_2 t}
\int_0^{\frac{t}{2}}e^{\alpha_2\tau}\left[C(1+\tau)^{-
\frac{2\theta'}{1-2\theta'}}+\|\nabla \Delta
\mathbf{d}_P(\tau)\|^2\right] d
\tau\non\\
&& +e^{-\alpha_2 t} \int_{\frac{t}{2}}^t
e^{\alpha_2\tau}\left[C(1+\tau)^{-
\frac{2\theta'}{1-2\theta'}}+\|\nabla \Delta
\mathbf{d}_P(\tau)\|^2\right] d \tau
\non\\
&:=& Ce^{-\alpha_2 t} + Z_1(t)+Z_2(t).\label{cwscr12}
 \eea
It follows from \eqref{A6} and (H6) that
 \bea Z_1(t)&\leq& Ce^{-\frac{\alpha_2}{2} t}
\int_0^{\frac{t}{2}}\left[C(1+\tau)^{-
\frac{2\theta'}{1-2\theta'}}+\|\nabla \Delta
\mathbf{d}_P(\tau)\|^2\right] d \tau\non\\
&\leq& Ce^{-\frac{\alpha_2}{2}
t}\left(t+\int_0^\frac{t}{2}(1+\tau)^{-2-2\gamma}d\tau\right)\non\\
&\leq &
C(1+t)^{- \frac{2\theta'}{1-2\theta'}}.\non
 \eea
 Next, by \eqref{A9} and the fact $\frac{2\theta'}{1-2\theta'}< 1+2\gamma $, we deduce that
 \bea
 Z_2(t)&\leq& Ce^{-\alpha_2 t} \left(1+\frac{t}{2}\right)^{-
\frac{2\theta'}{1-2\theta'}} \int_{\frac{t}{2}}^t
e^{\alpha_2\tau}d\tau +  C\int_{\frac{t}{2}}^t \|\nabla \Delta
\mathbf{d}_P(\tau)\|^2 d \tau\non\\
&\leq&  C(1+t)^{- \frac{2\theta'}{1-2\theta'}}+
C(1+t)^{-1-2\gamma}\non\\
&\leq& C(1+t)^{- \frac{2\theta'}{1-2\theta'}}.\non
 \eea
As a result, we obtain that
 \be
 z(t)\leq C(1+t)^{- \frac{2\theta'}{1-2\theta'}}, \quad \forall\
 t\geq 0.\label{Z}
 \ee
In particular, we have
  \be
  \mathcal{A}_P(t) \leq   C(1+t)^{- \frac{2\theta'}{1-2\theta'}},\quad
  t\geq 0,\label{rateA}
  \ee
  which together with  \eqref{eD} and \eqref{rate4} yields the
  following estimate
  \be
  \|\mathbf{v}(t)\|^2_{\mathbf{V}}+\|\Delta \mathbf{d}(t)-\Delta
  \mathbf{d}_\infty\|^2\leq   C(1+t)^{- \frac{2\theta'}{1-2\theta'}},\quad
  \forall\,t\geq 0.\non
  \ee
  Finally, using a standard elliptic estimate, we obtain (cf. (H7))
  \be
 \|\mathbf{d}(t)-
  \mathbf{d}_\infty\|_{\mathbf{H}^2}
  \leq C \|\Delta \mathbf{d}(t)-\Delta
  \mathbf{d}_\infty\|+
  C\|\mathbf{h}(t)-\mathbf{h}_\infty\|_{\mathbf{H}^\frac32(\Gamma)}
  \leq C(1+t)^{- \frac{\theta'}{1-2\theta'}},\non
  \ee
for all $t\geq 0$ and this finishes the proof.
\end{proof}


\section{Long-time behavior in 3D}
\setcounter{equation}{0}\label{Sec6}

As in the classical Navier--Stokes case (see \cite{Ler}), we can prove the eventual
regularity of any global weak solution. Thus the convergence results can also
be extended to the 3D case. Indeed, comparing with Lemma \ref{H3D}, we derive
first an alternative higher-order energy inequality.

\bl \label{H3Ds} Let the assumptions of Proposition
\ref{we} hold for all $T>0$. Suppose, in addition, that
\eqref{hyp6bis}--\eqref{hyp8bis} are satisfied. If a weak solution
$(\mathbf{v}, \mathbf{d})$ is smooth enough then it fulfills the
following inequality
  \be  \frac{d}{dt}\mathcal{A}_P(t)+ \nu\|S \mathbf{v}\|^2
  +\|\nabla (\Delta
  \wt{\mathbf{d}}-\mathbf{f}(\mathbf{d}))\|^2\leq
  C_*(\mathcal{A}^3_P(t)+\mathcal{A}_P(t) + R_3(t)),\label{h3dss}
  \ee
  where
  \be
  R_3(t)=  \|\p_t \mathbf{d}_P(t)\|^6+\|\p_t \mathbf{d}_P(t)\|^2
 +\|\nabla \Delta \mathbf{d}_P(t)\|^2+\|\mathbf{g}(t)\|^2,
  \ee
  for all $t\geq 0$. Here $C_*$ is a positive constant that may depend on
  $\nu$, $\|\mathbf{v}_0\|$, $\|\mathbf{d}_0\|_{\mathbf{H}^1}$,
  $\|\mathbf{g}\|_{L^2(0, +\infty; \mathbf{V}^*)}$,
  $ \|\mathbf{h}\|_{L^2_{tb}(0,+\infty;\mathbf{H}^\frac32(\Gamma))}$,
  $\|\mathbf{h}_t\|_{ L^2(0,+\infty; \mathbf{H}^{-\frac12}(\Gamma))}$, $\|\mathbf{h}_t\|_{L^1(0,+\infty;
 \mathbf{H}^{-\frac12}(\Gamma))}$.
 \el
 \begin{proof}
 We reconsider the estimates in the proof of Lemma \ref{H3D}.  Recalling \eqref{2dI1}
 and \eqref{split}, thanks to the Young inequality, it is not difficult to obtain that
 \be
 |I_1|\leq \|S \mathbf{v}\|\|\mathbf{v}\|_{\mathbf{L}^6}\|\nabla  \mathbf{v}\|_{\mathbf{L}^3}
 \leq C\|S \mathbf{v}\|^\frac32\|\nabla \mathbf{v}\|^\frac32\leq \varepsilon\|S \mathbf{v}\|^2+C\|\nabla \mathbf{v}\|^6,
 \non
 \ee
 \be
 |I_2|\leq
   \varepsilon \|S
   \mathbf{v}\|^2+C\|\mathbf{g}\|^2,\non
   \ee
   \bea
 |I_3|
  &\leq& \varepsilon\|S \mathbf{v}\|^2+\varepsilon\|\nabla (\Delta
  \wt{\mathbf{d}}-\mathbf{f}(\mathbf{d}))\|^2+C\|\Delta
  \wt{\mathbf{d}}-\mathbf{f}(\mathbf{d})\|^6+C\|\Delta
  \wt{\mathbf{d}}-\mathbf{f}(\mathbf{d})\|^2\non\\
  && +C\|\p_t\mathbf{d}_P\|^6+C\|\p_t\mathbf{d}_P\|^2+C\|\nabla
  \Delta\mathbf{d}_P\|^2,
  \non
  \eea
  \bea
 |I_4|
 &\leq &\varepsilon\|S \mathbf{v}\|^2+\varepsilon\|\nabla (\Delta
  \wt{\mathbf{d}}-\mathbf{f}(\mathbf{d}))\|^2+C\|\Delta
  \wt{\mathbf{d}}-\mathbf{f}(\mathbf{d})\|^6\non\\
  && +C\|\nabla \mathbf{v}\|^6+C\|\nabla \mathbf{v}\|^2 +C\|\p_t\mathbf{d}_P\|^6,\non
  \eea
  \be
 |I_{5a}| \leq \varepsilon\|S \mathbf{v}\|^2+C\|\Delta
  \wt{\mathbf{d}}-\mathbf{f}(\mathbf{d})\|^2+C\|\nabla
  \mathbf{v}\|^2+C\|\p_t\mathbf{d}_P\|^2.\non
 \ee
 In addition, $I_{5b}$ can be exactly estimated as \eqref{2dI5b}.
 Collecting all the estimates, and taking $\varepsilon$ to be sufficiently small,
 we obtain  our conclusion \eqref{h3dss}.
 \end{proof}

 Then we prove the following sufficient condition for the existence of global strong solution in 3D.

 \bp \label{aa3dloc}
Suppose that the assumptions of Proposition
\ref{we} and \eqref{hyp6bis}--\eqref{hyp8bis} are satisfied.
In addition, assume that $(\mathbf{v}_0, \mathbf{d}_0)\in \mathbf{V}\times
 \mathbf{H}^2(\Omega)$. If there exists a sufficiently small $\varepsilon_0\in (0,1]$ such that
 \be
 \int_0^{+\infty}(\nu\|\nabla \mathbf{v}(t)\|^2+\|\Delta \wha{
\mathbf{d}}(t)-\mathbf{f}(\mathbf{d}(t))\|^2)dt \leq
\varepsilon_0.\label{eqenergy}
 \ee
then problem \eqref{1}--\eqref{5} admits a unique global strong solution
$(\mathbf{v},\mathbf{d})$ in $\Omega\times(0,+\infty)$, provided
that $\|\mathbf{h}_t\|_{L^2(0,+\infty;\mathbf{H}^\frac12(\Gamma))}$ is small
enough.
 \ep
 \begin{proof}
 For simplicity, we give a formal proof. To make it rigorous we should work within a proper
approximation scheme (see, for instance, \cite{B,C09}).
 Let $L_i>0$ $(i=1,2,3,4,5)$ be the constants such that
 \bea
 \|\mathbf{v}_0\|+\|\mathbf{d}_0\|_{\mathbf{H}^1}&\leq& L_1,\label{loc1}\\
 \|\mathbf{h}_t\|_{L^2(0,+\infty;\mathbf{H}^\frac12(\Gamma))} &\leq& L_2,
 \label{loc2}\\
 \|\mathbf{h} \|_{L^2_{tb}(0,+\infty;\mathbf{H}^\frac32(\Gamma))}&\leq&L_3, \label{loc3}\\
 \|\mathbf{h}_t\|_{ L^1(0,+\infty; \mathbf{H}^{-\frac12}(\Gamma))}&\leq& L_4,\label{loc4}\\
 \|\mathbf{g}\|_{L^2(0, +\infty; \mathbf{V}^*)}&\leq&L_5.\label{loc5}
 \eea
  It follows from the basic energy inequality \eqref{EN1} that
 \be
 \wha{\mathcal{E}}(t)
 + \frac12\int_0^t \left( \nu \|\nabla \mathbf{v}\|^2+\|\Delta \wha{\mathbf{d}}-\mathbf{f}(\mathbf{d})\|^2\right)d\tau
 \leq \wha{\mathcal{E}}(0)+ \int_0^{+\infty} r(t) dt,\quad \forall\, t\geq 0.\label{eeee}
 \ee
 Then, by definition of $\wha{\mathcal{E}}$ and Lemma \ref{Ap1}, we have
   \bea
 && \|\mathbf{v}(t)\|+\|\mathbf{d}(t)\|_{\mathbf{H}^1}\leq C_1,\quad
 \forall\, t\geq 0,\label{kkk1}\\
 && \int_0^{+\infty} \left(\nu \|\nabla \mathbf{v}\|^2+\|\Delta \wha{\mathbf{d}}-\mathbf{f}(\mathbf{d})\|^2\right)dt\leq C_2, \label{kkkk}
 \eea
 where the constants $C_1, C_2$ depend on $L_1,...,L_5$ and $\Omega$.

Let $K>0$ be such that
 \be
 \nu\|\nabla \mathbf{v}_0\|^2+\|\Delta
\wt{\mathbf{d}}(0)-f(\mathbf{d}_0)\|^2\leq K.\label{loc6}
 \ee
  Keeping Lemma \ref{H3Ds} in mind and arguing as in
 \cite{LL95}, we consider the following Cauchy problem
 \be
  \frac{d}{dt}Y(t)=C_*(Y(t)^3+Y(t))+C_*R_3(t),\quad
Y(0)= \max\left\{1,\nu^{-1}\right\}K\geq
\mathcal{A}_P(0).\label{ode}
 \ee
 We denote by $I=[0,T_{max})$ the (right) maximal interval for the existence
 of a (nonnegative) solution $Y(t)$ so that
 $ \displaystyle\lim_{t\rightarrow T_{max}^-} Y(t)=+\infty.$
 It easily follows from \eqref{h3dss} and the comparison principle
 that $0\leq\mathcal{A}_P(t)\leq Y(t)$, for any $t\in I$.
 Consequently, $\mathcal{A}_P(t)$ is finite on $I$.
  We deduce from Lemma \ref{Ap2} that
 \be
 \int_0^{+\infty} R_3(t) dt\leq C_3,\non
 \ee
 where $C_3$ is a constant depending on $\Omega$,
  $\|\mathbf{g}\|_{L^2(0, +\infty; \mathbf{H})}$ and
  $L_2$.
 Besides, we note that $T_{max}$ is determined by $Y(0)$, $C_*$ and $C_3$ such that $T_{max}=T_{max}(Y(0),C_*, C_3)$ is
 increasing when $Y(0)\geq 0$ is decreasing.
  Taking $t_0=\frac12 T_{max}> 0$, then it
 follows that $Y(t)$ (as well as $\mathcal{A}_P(t)$) is uniformly bounded on
$[0, t_0]$. This easily implies the local existence of a unique strong
solution to problem \eqref{1}--\eqref{5} (at least) on $[0,t_0]$
(actually on $[0, T_{max})$, but we lose uniform estimates on such
maximal interval).

By Lemma \ref{Ap2} (cf. \eqref{A4}), we have
 \be
 \sup_{t\geq 0} \|\Delta
 (\mathbf{d}_P(t)-\mathbf{d}_E(t))\|^2\leq c\|\mathbf{h}_t\|_{L^2(0,+\infty;\mathbf{H}^\frac12(\Gamma))}^2,\label{A4aa}
 \ee
 where $c$ is a constant that depends only on $\Omega$.
Set now
 \be
 \bar{\varepsilon}_0= \min\left\{1, \frac{t_0K}{8}\right\}, \quad
 L_6=\min\left\{1,L_2^2,
 \frac{K}{4c}\right\}.
 \label{choice1}
 \ee
 From the assumption, there exists a small constant $\varepsilon_0\leq \bar{\varepsilon}_0$ such that \eqref{eqenergy} is satisfied. Therefore, we can find $t_* \in [\frac{t_0}{2}, t_0]$ such that
 \be
  \nu\|\nabla \mathbf{v}(t_*)\|^2+\|\Delta \wha{\mathbf{d}}(t_*)-\mathbf{f}(\mathbf{d}(t_*))\|^2
\leq 2\bar{\varepsilon}_0t_0^{-1}. \non
 \ee
 Moreover, if we further assume $$\|\mathbf{h}_t\|_{L^2(0,+\infty;\mathbf{H}^\frac12(\Gamma))}^2\leq L_6,$$
 then by \eqref{A4aa}  we obtain
  \bea
  \mathcal{A}_P(t_*)&\leq & \nu\|\nabla \mathbf{v}(t_*)\|^2 +\|\Delta
 \wt{\mathbf{d}}(t_*)-\mathbf{f}(\mathbf{d}(t_*))\|^2\non\\
 &\leq & \nu\|\nabla \mathbf{v}(t_*)\|^2+ 2\|\Delta
 \wha{\mathbf{d}}(t_*)-\mathbf{f}(\mathbf{d}(t_*))\|^2+ 2\|\Delta
 (\mathbf{d}_P(t_*)-\mathbf{d}_E(t_*))\|^2\non\\
 &\leq& \nu\|\nabla \mathbf{v}(t_*)\|^2+2\|\Delta
 \wha{\mathbf{d}}(t_*)-\mathbf{f}(\mathbf{d}(t_*))\|^2
 +2c\|\mathbf{h}_t\|_{L^2(0,+\infty;\mathbf{H}^\frac12(\Gamma))}^2
 \non\\
 &\leq& 4\bar{\varepsilon}_0t_0^{-1}+2cL_6\leq K\non\\
 &\leq&
 \max\left\{1,\nu^{-1}\right\}K=Y(0).\non
 \eea
 Taking $t_*$ as the initial time for the ordinary differential equation \eqref{ode}, we
infer from the above argument that
$\mathcal{A}_P(t)$ is uniformly bounded at least on $[0,
\frac{3t_0}{2}]\subset [0,t_*+t_0]$. Moreover, its bound only
depends on $\Omega$, $\nu$, $L_1,...,L_6$, $C_*$ and
$t_0$. Then by an iterative argument we can show that $\mathcal{A}_P(t)$ is
uniformly bounded for all $t\geq 0$ and this enable us to extend the
local strong solution to the whole time interval $[0,+\infty)$. The proof is
complete.
 \end{proof}

A consequence of the above proposition is the eventual regularity of global weak solutions.
\bt \label{reg}
Suppose that the assumptions of Proposition
\ref{we} and \eqref{hyp6bis}--\eqref{hyp8bis} are satisfied.
 Let $(\mathbf{v}, \mathbf{d})$ be a global weak solution to \eqref{1}--\eqref{5}.
 Then there exists a large time $T^*\in (0, +\infty)$ such that  $(\mathbf{v}, \mathbf{d})$ is a strong solution
 on $(T^*,+\infty)$.
 \et
 \begin{proof}
 Let $L_1, L_2, L_3, L_4, L_5>0$ be the constants as in the proof of Proposition \ref{a3dloc}. For a weak solution $(\mathbf{v}, \mathbf{d})$, we still have the uniform estimates \eqref{kkk1} and \eqref{kkkk}.
Considering the ODE problem \eqref{ode}, we can fix the constants $\bar{\varepsilon}_0$, $L_6$ and $t_0$.
Taking $\varepsilon_0=\bar{\varepsilon}_0$, we observe that
  there must exist a sufficiently large $T_1>0$ such that
  \bea
   \int_{T_1}^{+\infty} \left(\nu \|\nabla \mathbf{v}\|^2+\|\Delta \wha{\mathbf{d}}-\mathbf{f}(\mathbf{d})\|^2\right)dt &\leq& \varepsilon_0,\label{t1}\\
 \|\Delta \mathbf{d}_P(t)-\Delta \mathbf{d}_E(t)\|&\leq& L_6,\quad \forall\, t\geq [T_1,+\infty),\label{l6a}
    \eea
    where for the second inequality we have used Lemma \ref{Ap2}(i) and the fact that $\partial_t \mathbf{d}_P(t)=\Delta \mathbf{d}_P(t)-\Delta \mathbf{d}_E(t)$.
    Also, \eqref{t1} implies that there is $T^*\in [T_1, T_1+2t_0]$ such that
    \be
     \nu\|\nabla \mathbf{v}(T_*)\|^2+ \|\Delta
 \wha{\mathbf{d}}(T_*)-\mathbf{f}(\mathbf{d}(T_*))\|^2\leq \frac{\bar{\varepsilon}_0}{t_0}.
    \ee
    As a result,
    \bea
  && \nu\|\nabla \mathbf{v}(T_*)\|^2 +\|\Delta
 \wt{\mathbf{d}}(T_*)-\mathbf{f}(\mathbf{d}(T_*))\|^2\non\\
 &\leq & \nu\|\nabla \mathbf{v}(T_*)\|^2+ 2\|\Delta
 \wha{\mathbf{d}}(T_*)-\mathbf{f}(\mathbf{d}(T_*))\|^2+ 2\|\Delta
 (\mathbf{d}_P(T_*)-\mathbf{d}_E(T_*))\|^2\non\\
  &\leq& \frac{2\bar{\varepsilon}_0}{t_0}+2cL_6\non\\
 &\leq&
  K.\non
 \eea
  Taking $T^*$ as the initial time, then we can apply Proposition \ref{aa3dloc} to conclude that problem  \eqref{1}--\eqref{5} admits a unique global strong solution
  $(\mathbf{v}', \mathbf{d}')$. By the weak/strong uniqueness result \cite[Theorem 7]{C09}, we see that $(\mathbf{v}, \mathbf{d})$ coincides with $(\mathbf{v}', \mathbf{d}')$
  on $[T^*, +\infty)$. The proof is complete.
\end{proof}

Thanks to the eventual regularity result we can argue as in the previous section to prove the following result

 \bt \label{conweak}
  Suppose that the assumptions of Theorem \ref{reg} hold.
  Then any global weak solution given by Proposition \ref{we}
  converges in $\mathbf{V}\times \mathbf{H}^2(\Omega)$
  to a single equilibrium $(\mathbf{0},\mathbf{d}_\infty)$ with estimates on the convergence rate similar to the 2D case, provided that $\mathbf{g}$ and $\mathbf{h}$ fulfill the corresponding hypotheses (H1)--(H7) as in Theorems \ref{conv1} and \ref{convrate1}.
  \et

\br
We recall that there exists a (unique) global strong solution when the viscosity is large enough (cf. Theorem \ref{exe3d}).
Consequently, due to Lemma \ref{H3D}, all the results proven in Section \ref{Sec5} (i.e., Theorem \ref{conv1} and
Theorem \ref{convrate1}) still hold with the same assumptions on the data. The related proofs just require some minor modifications.
\er

The existence of a global strong solution is also ensured (with no restrictions on the fluid viscosity)
when the initial data are close to a given equilibrium and
the time dependent boundary data satisfies suitable bounds. First, recall that the basic energy
inequality \eqref{EN1} implies  (cf. \eqref{eeee})
 \be
\int_0^{t}(\nu\|\nabla \mathbf{v}(t)\|^2+\|\Delta \wha{
\mathbf{d}}(t)-\mathbf{f}(\mathbf{d}(t))\|^2)dt \leq  2(\wha{\mathcal{E}}(0)- \wha{\mathcal{E}}(t))+2\int^{+\infty}_0r(t) dt,\non
 \ee
 and
 \be\int^{+\infty}_0r(t) dt\leq C_r\left( \|\mathbf{h}_t\|^2_{L^2(0,+\infty;\mathbf{H}^{-\frac12}(\Gamma))}+\|\mathbf{h}_t\|_{ L^1(0,+\infty; \mathbf{H}^{-\frac12}(\Gamma))}+\|\mathbf{g}\|^2_{L^2(0, +\infty; \mathbf{V}^*)}\right),\label{esonr}
 \ee
 where $C_r$ is a universal constant. Then we can easily deduce from Proposition \ref{aa3dloc}
that if the
lifted energy stays sufficiently close to its initial state, then system \eqref{1}--\eqref{5} admits a unique global strong solution (cf.
\cite{LL95} for the autonomous case).

 \bp \label{a3dloc}
Assume \eqref{hyp6bis}--\eqref{hyp8bis} and \eqref{hyp4} hold.
Moreover, suppose that $(\mathbf{v}_0, \mathbf{d}_0)\in
\mathbf{V}\times
 \mathbf{H}^2(\Omega)$ satisfying \eqref{hyp5} and
 $|\mathbf{d}_0|_{\mathbb{R}^3}\leq 1$.
If there exists a sufficiently small $\varepsilon_0\in (0,1]$ such that
 \be
\wha{\mathcal{E}}(t)\geq \wha{\mathcal{E}}(0)-\varepsilon_0,
\quad \forall\, t\geq 0, \label{eqenergy1}
 \ee where
$\wha{\mathcal{E}}$ is the lifted energy defined by \eqref{E},
then problem \eqref{1}--\eqref{5} admits a unique global strong solution
$(\mathbf{v},\mathbf{d})$ in $\Omega\times(0,+\infty)$, provided
that $\|\mathbf{h}_t\|_{L^2(0,+\infty;\mathbf{H}^\frac12(\Gamma))}$, $\|\mathbf{h}_t\|_{ L^1(0,+\infty; \mathbf{H}^{-\frac12}(\Gamma))}$ and $\|\mathbf{g}\|_{L^2(0, +\infty; \mathbf{V}^*)}$ are small
enough.
 \ep

Let us assume that for all $t\geq 0$ (comparing with assumptions (H1), (H4), (H5))
\begin{itemize}
\item[(H1')]
    $\int_t^{+\infty}\|\mathbf{h}_t(\tau)\|_{\mathbf{H}^\frac12(\Gamma)}d\tau\leq
    M_1(1+t)^{-1-\gamma}$;

\item[(H4')] $ \|\mathbf{g}(t)\|^2\leq
    M_2(1+t)^{-2-\gamma}$;

 \item[(H5')]
     $\|\mathbf{h}_t(t)\|_{\mathbf{L}^2(\Gamma)}\leq
     M_3(1+t)^{-1-\gamma}$.
\end{itemize}
 Here $M_j$, $j=1,2,3$ and $\gamma$ are positive
constants. $\gamma$ characterizes the decay rate of
non-autonomous terms, while $M_j$ control their magnitude.

In spirit of Proposition \ref{a3dloc}, in what follows, we prove the global existence of a strong solution that
originates near a local minimizer of the lifted energy with suitably small perturbations in terms of the nonautonomous terms $\mathbf{h}$ and $\mathbf{g}$ (namely, the magnitudes $M_j$ should be sufficiently small).

\bt \label{3dlom} Suppose that \eqref{hyp6bis}--\eqref{hyp8bis} and \eqref{hyp4} hold, the constant
$\gamma>1$.
Moreover, assume that $(\mathbf{v}_0, \mathbf{d}_0)\in
\mathbf{V}\times
 \mathbf{H}^2(\Omega)$ satisfies \eqref{hyp5} and
 $|\mathbf{d}_0|_{\mathbb{R}^3}\leq 1$. Denote by
$\mathbf{d}_{E}^*$ the unique solution to
 \be
 \begin{cases}
 -\Delta \mathbf{d}^*_E=\mathbf{0},\quad x\in \Omega,\\
 \mathbf{d}_E^*=\mathbf{h}_\infty, \quad x\in \Gamma,
 \end{cases}\label{dEs}
 \ee
and set
 $$\mathscr{E}(\mathbf{d})=\frac12\|\nabla
 (\mathbf{d}-\mathbf{d}^*_E)\|^2+\int_\Omega F(\mathbf{d})dx, \quad \forall \mathbf{d}\in \mathcal{N}.$$
  Let $\mathbf{d}^*\in \mathcal{N}\cap\,\mathbf{H}^2(\Omega)$ be a local minimizer of
 $\mathscr{E}(\mathbf{d})$
 in the sense that $\mathscr{E}(\mathbf{d})
 \geq \mathscr{E}(\mathbf{d}^*)$ for all $\mathbf{d}\in \mathcal{N}$
 satisfying
 $\|\mathbf{d}-\mathbf{d}^*\|_{\mathbf{H}^1}<\delta$, where $\delta>0$ is a certain small constant.
 Suppose also that the initial data $\mathbf{v}_0$ and $\mathbf{d}_0$
satisfy
 \be
 \|\mathbf{v}_0\|_{\mathbf{V}}\leq 1, \quad
 \|\mathbf{d}_0-\mathbf{d}^*\|_{\mathbf{H}^2}\leq 1.
 \label{localbd1}
 \ee
Then there exist positive constants
 $\sigma_1, \sigma_2, M_1, M_2, M_3, L_0$, which are sufficiently small and may depend
 on the system coefficients, on $\Omega$ and on $\mathbf{d}^*$,
such that if the initial data $(\mathbf{v}_0,\mathbf{d}_0)$ and $\mathbf{h}$ also
fulfill
 \be
 \|\mathbf{v}_0\|\leq \sigma_1, \quad
 \|\mathbf{d}_0-\mathbf{d}^*\|_{\mathbf{H}^1}\leq \sigma_2,\quad
 \|\mathbf{h}_t\|_{L^2(0,+\infty;\mathbf{H}^\frac12(\Gamma))}^2 \leq
 L_0, \non
 \ee
  and (H1'), (H4'), (H5') hold with such $M_j$, $j=1,2,3$, then problem
 \eqref{1}--\eqref{5} admits a unique global strong solution $(\mathbf{v}, \mathbf{d})$.
 \et
 \begin{proof}
 Without loss of generality, we assume $\delta\in(0,1]$.
 In the subsequent proof, $C_i$ ($i\in\mathbb{N}$) stand for positive constants that only depend
 on $\Omega$, $\nu$, $\gamma$ and $\mathbf{d}^*$.
 Under the current assumption \eqref{localbd1} on the initial data,  it is not difficult to
 see that the constants $L_1$ and $K$ in \eqref{loc1} and \eqref{loc6}
 depend on $\mathbf{d}^*$ only. We just take $L_2=L_3=L_4=L_5=1$ in \eqref{loc3}
 for the sake of simplicity. Then we have
 the uniform estimate (cf. \eqref{kkk1})
 \be
 \|\mathbf{v}(t)\|+\|\mathbf{d}(t)\|_{\mathbf{H}^1}\leq C_1,\quad
 t\geq 0.\non
 \ee
Arguing as in the proof of Proposition \ref{aa3dloc}, we find that
problem \eqref{1}--\eqref{5}
admits a unique strong solution (at least) on $[0,t_0]$, whose
$\mathbf{V}\times\mathbf{H}^2$ norm is uniformly bounded on
$[0,t_0]$:
 \be
 \|\mathbf{v}(t)\|_{\mathbf{V}}+\|\mathbf{d}(t)\|_{\mathbf{H}^2}\leq C_3,\quad
 t\in[0,t_0]. \label{lochigh}
 \ee
Besides, we can also fix the constants $\bar{\varepsilon}_0$ and
 $L_6$ (see \eqref{choice1}). In the subsequent proof, we just take $$\varepsilon_0=\bar{\varepsilon}_0,\quad  L_0=L_6.$$
It follows from \eqref{esonr} that
$$\int_0^{+\infty}r(t)dt\leq C_rC_s(M_1+M_2+M_3^2)\leq \frac{\varepsilon_0}{4},$$
provided that $M_1, M_2, M_3>0$ are assumed to be properly small and satisfying
$$ M_1+M_2+M_3^2\leq \frac{\varepsilon_0}{4C_rC_s},$$
where $C_s$ is a universal constant due to the Sobolev embedding.
Hence, according to Propositions \ref{aa3dloc} and \ref{a3dloc}, in order to prove the existence of global strong
solution, we only have to
 verify that
 \be
 \wha{\mathcal{E}}(t)-\wha{\mathcal{E}}(0)\geq
 -\frac{\varepsilon_0}{2}, \quad \forall\, t\geq 0.\label{egoal}
 \ee
 First, we notice that (recalling  \eqref{LE}, \eqref{iE0} and \eqref{fEa})
 \bea
 && \wha{\mathcal{E}}(0)-\wha{\mathcal{E}}(t)\non\\
 &\leq&
 \frac12\|\mathbf{v}_0\|^2+
 \wha{E}(\mathbf{d}_0)-\wha{E}(\mathbf{d}(t))\non\\
 &=& \frac12\|\mathbf{v}_0\|^2+ \frac12\|\nabla
 (\mathbf{d}_0-\mathbf{d}_{E0})\|^2-\frac12\|\nabla
 (\mathbf{d}(t)-\mathbf{d}_{E})\|^2+\int_\Omega
 F(\mathbf{d}_0)-F(\mathbf{d}(t)) dx\non\\
 &\leq& \frac12\|\mathbf{v}_0\|^2+
 C_4(\|\mathbf{d}_0-\mathbf{d}(t)\|_{\mathbf{H}^1}+\|\mathbf{d}_{E0}-\mathbf{d}_E\|_{\mathbf{H}^1}).\label{a1}
 \eea
 On the other hand, thanks to standard elliptic estimates, we have
 \bea
 \|\mathbf{d}_{E0}-\mathbf{d}_E\|_{\mathbf{H}^1}&\leq&
 c\|\mathbf{d}_0|_\Gamma-\mathbf{h}(t)\|_{\mathbf{H}^\frac12(\Gamma)}\non\\
 &\leq&
 c\|\mathbf{d}_0|_\Gamma-\mathbf{h}_\infty\|_{\mathbf{H}^\frac12(\Gamma)}
 +c\|\mathbf{h}_\infty-\mathbf{h}(t)\|_{\mathbf{H}^\frac12(\Gamma)}\non\\
 &\leq& c\|\mathbf{d}_0-\mathbf{d}^*\|_{\mathbf{H}^1}+
 c\int_t^{+\infty}\|\mathbf{h}_t(\tau)\|_{\mathbf{H}^\frac12(\Gamma)}d\tau\non\\
 &\leq& c\sigma_2+cM_1, \quad \forall\ t\geq 0.\label{a2}
 \eea
 Let
 \be
 \sigma_1\leq \min\left\{1, \frac{\sqrt{\varepsilon_0}}{2}\right\},\quad
 \sigma_2\leq \frac{\varepsilon_0}{8C_4}\min\{1,c^{-1}\}, \quad M_1\leq
 \min\left\{1,\frac{\varepsilon_0}{8C_4c}\right\}.\label{m1}
 \ee
 Due to \eqref{a1} and \eqref{a2},  in order to prove \eqref{egoal}, we only have to verify
 \be
 \|\mathbf{d}_0-\mathbf{d}(t)\|_{\mathbf{H}^1}\leq
 \frac{\varepsilon_0}{8C_4},\quad \forall\ t\geq 0. \label{egoal1}
 \ee
 Since $\mathbf{d}^*\in \mathcal{N}\cap\, \mathbf{H}^2(\Omega)$ is the local minimizer of $\mathscr{E}$, it is
 easily to verify that $\mathbf{d}^*$ satisfies \eqref{sta} and thus is the critical point of $E$.
 As a consequence, Corollary \ref{ELS} holds for $\mathbf{d}^*$ with
 constants
 $\theta, \beta$ determined by $\mathbf{d}^*$. By \eqref{theta},
 $\theta'$ can be determined by $\theta$ and $\gamma$. In addition,
 we further choose $\theta'$ smaller if necessary such that (recall
 that $\gamma >1$)
 \be
 \theta'\leq \frac{\gamma-1}{2\gamma}.\label{theta1}
 \ee
Let us define
 \be
 \varpi=\min\left\{\frac{\beta}{2},\  \frac{\delta}{2},\
 \frac{\varepsilon_0}{10C_4}\right\},\label{vpi}
 \ee
and set
 \be
 \bar t_0=\sup\{t\in [0,t_0], \
 \|\mathbf{d}(t)-\mathbf{d}^*\|_{\mathbf{H}^1}<\varpi,\  \forall\
 s\in [0,t)\}\non
 \ee
 If we assume
 \be
 \sigma_2\leq \frac14\varpi,
 \ee
 then by the continuity of $\mathbf{d}(t)$ in $\mathbf{H}^1(\Omega)$, we have $\bar
 t_0>0$. Next, we shall prove that $\bar t_0>t_0$ by contradiction.
We introduce the auxiliary functional
 \be
 \Psi_1(t)=\wha{\mathcal{E}}(t)-\wha{E}(\mathbf{d}^*)+2\int_t^{+\infty}r(\tau)d\tau,\non
 \ee
 and the function
 \be
 \bar {\mathbf{d}}(t)=\mathbf{d}(t)-\mathbf{d}_E+\mathbf{d}_E^*.\non
 \ee
 It easily follows that
 \bea
 \Psi_1(t)&\geq& \wha{E}(\mathbf{d}(t))-\wha{E}(\mathbf{d}^*)=
  \wha{E}(\mathbf{d}(t))-\mathscr{E}(\bar {\mathbf{d}}(t))+\mathscr{E}(\bar {\mathbf{d}}(t))-\wha{E}(\mathbf{d}^*)\non\\
  &=& \int_\Omega F(\mathbf{d}(t))-
 F(\bar {\mathbf{d}}(t))dx+\mathscr{E}(\bar
 {\mathbf{d}}(t))-\wha{E}(\mathbf{d}^*).\label{Psia}
 \eea
 By definition, $\bar {\mathbf{d}}(t)\in \mathcal{N}$. Moreover, on
 $[0,\bar t_0]$,
 \bea
 && \|\bar {\mathbf{d}}(t)-\mathbf{d}^*\|_{\mathbf{H}^1}\leq
 \|\mathbf{d}(t)-\mathbf{d}^*\|_{\mathbf{H}^1}+\|\mathbf{d}_E-\mathbf{d}_E^*\|_{\mathbf{H}^1}\non\\
 &\leq&
 \varpi+c\|\mathbf{h}(t)-\mathbf{h}_\infty\|_{\mathbf{H}^\frac12(\Gamma)}\leq
 \frac{\delta}{2}+c\int_t^{+\infty}\|\mathbf{h}_t(\tau)\|_{\mathbf{H}^\frac12(\Gamma)}d\tau\non\\
 &\leq& \frac{\delta}{2}+cM_1.\non
 \eea
 Taking
 \be
 M_1\leq \min\left\{1,\frac{\delta}{4c}\right\},\label{M1}
 \ee
 then we have $\|\bar {\mathbf{d}}(t)-\mathbf{d}^*\|_{\mathbf{H}^1}\leq
 \delta$. Since $\mathbf{d}^*$ is a local minimizer of $\mathscr{E}$, we see that
 \be
 \mathscr{E}(\bar {\mathbf{d}}(t))- \wha{E}(\mathbf{d}^*)=\mathscr{E}(\bar {\mathbf{d}}(t))-\mathscr{E} (\mathbf{d}^*)\geq 0, \quad t\in
 [0,\bar t_0].\label{Psib}
 \ee
 On the other hand,  since $|\mathbf{d}(t)|_{\mathbb{R}^3}\leq 1 $ and $|\bar {\mathbf{d}}(t)|_{\mathbb{R}^3}\leq 3$ (this is due to the maximum principle \eqref{max}),
 we infer from the standard elliptic estimate and (H5') that
 \bea
 \left|\int_\Omega F(\mathbf{d}(t))-
 F(\bar {\mathbf{d}}(t))dx\right|&\leq&
 C_5\|-\mathbf{d}_E+\mathbf{d}_E^*\|\non\\
 &\leq&
 C_5c\int_t^{+\infty}\|\mathbf{h}_t(\tau)\|_{\mathbf{L}^2(\Gamma)}d\tau\non\\
 & \leq& C_5c M_3\gamma^{-1}(1+t)^{-\gamma}.\label{Psic}
  \eea
 Let us introduce now two further functions
 \bea
 z(t)= (C_5c+1) M_3\gamma^{-1}(1+t)^{-\gamma},\quad
 \Psi(t)=\Psi_1(t)+z(t).\non
 \eea
 We deduce from \eqref{Psia}--\eqref{Psic} that
 \be
 \Psi(t)\geq M_3\gamma^{-1}(1+t)^{-\gamma}>0, \quad t\in[0,\bar
 t_0],\non
 \ee
 and by the basic energy inequality \eqref{EN1}
 \bea
 \frac{d}{dt}\Psi(t)&=&\frac{d}{dt}\wha{\mathcal{E}}(t)-2r(t)-(C_5c+1)M_3(1+t)^{-1-\gamma}\non\\
 &\leq& -\frac{1}{4}\min\{\nu,
 1\}\mathcal{D}^2(t) -(C_5c+1)M_3(1+t)^{-1-\gamma}\non\\
 &\leq& -
 C_6\left(\mathcal{D}(t)+M_3^\frac12(1+t)^{-\frac{1+\gamma}{2}}\right)^2,\non
 \eea
 where $ \mathcal{D}(t)=\|\nabla
 \mathbf{v}(t)\|+\|\Delta
 \wha{\mathbf{d}}(t)-\mathbf{f}(\mathbf{d}(t))\|$.
Arguing as to get \eqref{dePhi}, using Remark \ref{st} and
assumptions (H1'), (H4'), we deduce
 \be
 \Psi(t)^{1-\theta'}\leq C_7\left(\mathcal{D}(t)+ (M_1+M_2)
(1+t)^{-(1-\theta')(1+\gamma)}+
M_3^{1-\theta'}(1+t)^{-(1-\theta')\gamma}\right).\non
 \ee
 Assuming
 \be
 M_1\leq \frac12 M_3^\frac12, \quad M_2\leq \frac12 M_3^\frac12, \quad
 M_3\leq 1,
 \label{M123}
 \ee
we can see that
 \be
 \Psi(t)^{1-\theta'}\leq C_7\left(\mathcal{D}(t)+
 2M_3^\frac12(1+t)^{-(1-\theta')\gamma}\right).\non
 \ee
 As a result, we find
 \bea
 -\frac{d}{dt}\Psi(t)^{\theta'}&=&-\theta'
 \Psi(t)^{\theta'-1}\frac{d}{dt}\Psi(t)\non\\
 & \geq &
 \frac{C_6\left(\mathcal{D}(t)+M_3^\frac12(1+t)^{-\frac{1+\gamma}{2}}\right)^2}{C_7\left(\mathcal{D}(t)+M_3^\frac12
(1+t)^{-(1-\theta')\gamma}\right)}\non\\
 &\geq&
 C_8\left(\mathcal{D}(t)+M_3^\frac12(1+t)^{-\frac{1+\gamma}{2}}\right),\label{dPsi}
 \eea
 where we have used the fact that $\frac{1+\gamma}{2}\leq (1-\theta')\gamma$ (cf.  \eqref{theta1}).
 It follows from \eqref{vadt}, \eqref{lochigh},
 \eqref{M123}, \eqref{dPsi}, assumptions (H1'), (H4'), (H5') and the definition of $\Psi$ that
 \bea
  && \int_0^{\bar t_0} \|\mathbf{d}_t(t)\|dt \leq C_9
 \Psi(0)^{\theta'}\non\\
 & \leq&
 C_{10}\left(\|\mathbf{v}_0\|^{2}+\|\mathbf{d}_0-\mathbf{d}^*\|_{\mathbf{H}^1}
 +\|\mathbf{d}_{E0}-\mathbf{d}^*_E\|_{\mathbf{H}^1} +
 \int_0^{+\infty} r(t)dt+ z(0)\right)^{\theta'}\non\\
 &\leq&
 C_{11}\left(\|\mathbf{v}_0\|^{2}+\|\mathbf{d}_0-\mathbf{d}^*\|_{\mathbf{H}^1}+
 M_3^\frac12\right)^{\theta'}.\label{difff}
 \eea
 By \eqref{lochigh}, \eqref{difff} and an interpolation inequality, we get
 \bea
 && \|\mathbf{d}(\bar t_0)-\mathbf{d}^*\|_{\mathbf{H}^1}\non\\
  &\leq &
 \|\mathbf{d}(\bar
 t_0)-\mathbf{d}_0\|_{\mathbf{H}^1}+\|\mathbf{d}_0-\mathbf{d}^*\|_{\mathbf{H}^1}\non\\
 &\leq& C_{12}(\|\mathbf{d}(\bar
 t_0)\|_{\mathbf{H}^2}+ \|\mathbf{d}_0\|_{\mathbf{H}^2})^\frac12\|\mathbf{d}(\bar
 t_0)-\mathbf{d}_0\|^\frac12+\|\mathbf{d}_0-\mathbf{d}^*\|_{\mathbf{H}^1}\non\\
 &\leq&
 C_{13}\left(\|\mathbf{v}_0\|^{\theta'}+\|\mathbf{d}_0-\mathbf{d}^*\|_{\mathbf{H}^1}^{\frac{\theta'}{2}}+
 M_3^{\frac{\theta'}{4}}\right)+\|\mathbf{d}_0-\mathbf{d}^*\|_{\mathbf{H}^1}.\label{est0}
 \eea
 Taking now
 \bea
 && \sigma_1\leq \min\left\{1, \frac{\sqrt{\varepsilon_0}}{2},
 \left(\frac{\varpi}{6C_{13}}\right)^\frac{1}{\theta'}\right\},\quad
 \sigma_2\leq\min\left\{1,\frac14\varpi,
 \left(\frac{\varpi}{6C_{13}}\right)^\frac{2}{\theta'}\right\},\\
 && M_3\leq
 \min \left\{1, \left(\frac{\varpi}{6C_{13}}\right)^\frac{4}{\theta'}\right\},\label{m3}
 \eea
 we infer from \eqref{est0} that
 \be
 \|\mathbf{d}(\bar t_0)-\mathbf{d}^*\|_{\mathbf{H}^1}\leq
 \frac{3}{4}\varpi<\varpi.\non
 \ee
 This leads to a contradiction with the definition of $\bar t_0$.
 As
 a result, we have $\bar t_0> t_0$, and
 \bea
 \|\mathbf{d}_0-\mathbf{d}(t)\|_{\mathbf{H}^1} &\leq&
 \|\mathbf{d}_0-\mathbf{d}^*\|_{\mathbf{H}^1}+
 \|\mathbf{d}^*-\mathbf{d}(t)\|_{\mathbf{H}^1}\non\\
 &\leq&
 \sigma_2+\varpi\leq \frac54\varpi\leq \frac{\varepsilon_0}{8C_4},
 \quad \forall \ t\in [0,t_0].
 \eea
 Thus, \eqref{egoal} holds on $[0,t_0]$, which implies
  \be
\int_0^{t_0}(\nu\|\nabla \mathbf{v}(t)\|^2+\|\Delta \wha{
\mathbf{d}}(t)-\mathbf{f}(\mathbf{d}(t))\|^2)dt \leq
\varepsilon_0.\non
 \ee
  As in Proposition \ref{aa3dloc}, there exists $t_* \in [\frac{t_0}{2}, t_0]$ such that
 \be
  \nu\|\nabla \mathbf{v}(t_*)\|^2+\|\Delta \wha{\mathbf{d}}(t_*)-\mathbf{f}(\mathbf{d}(t_*))\|^2
\leq 2\varepsilon_0t_0^{-1},  \non
 \ee
 and again we have $
  \mathcal{A}_P(t_*)\leq \max\left\{1,\nu^{-1}\right\}K$.
 Taking $t_*$ as the initial time for the Cauchy problem \eqref{ode}, we can
 extend the (unique) strong solution to $[0,\frac32 t_0]$ and its $\mathbf{V}\times \mathbf{H}^2$-norm
  is uniformly bounded by the same constant $C_3$ as on $[0,t_0]$.
 Repeating the above argument in $[0,\frac32 t_0]$, we can verify that \eqref{egoal} still holds.
 By iteration we can show that \eqref{egoal} holds for all $t\geq
 0$. Hence, our conclusion follows from Proposition \ref{a3dloc}.
\end{proof}

Finally, we can conclude with the following local stability result:
  \bt \label{conloc}
  Let the assumptions of Theorem \ref{3dlom} hold.
  Then any global strong solution given by Theorem \ref{3dlom}
  converges in $\mathbf{V}\times \mathbf{H}^2(\Omega)$
  to a single equilibrium $(\mathbf{0},\mathbf{d}_\infty)$
  with $\mathbf{d}_\infty \in \mathcal{N}\cap\, \mathbf{H}^2(\Omega)$
  such that $\mathscr{E}(\mathbf{d}_\infty)=\mathscr{E}(\mathbf{d}^*)$.
  In addition, convergence rate estimates similar to the 2D case hold provided
  that $\mathbf{g}$ and $\mathbf{h}$ fulfill the corresponding hypotheses (i.e.,
  assumptions (H1), (H4), (H5) are replaced by (H1'), (H4'), (H5'), respectively).
  Indeed, the local energy minimizer $\mathbf{d}^*$ is (locally) Lyapunov stable, and in
particular, if $\mathbf{d}^*$ is an isolated local minimizer of
  $\mathscr{E}$, then it is (locally) asymptotically stable.
  \et
  \begin{proof}
 Arguing as in Section \ref{Sec5} we still find
 \be
  \displaystyle\lim_{t \to +\infty} \left(\|\mathbf{v}(t)\|_{\mathbf{V}}+
  \|\mathbf{d}(t)-\mathbf{d}_\infty\|_{\mathbf{H}^2}\right)=0,
  \label{conaaa}
 \ee
  for some $\mathbf{d}_\infty \in \mathcal{N}\cap\, \mathbf{H}^2(\Omega)$.
  The estimate on the convergence rates can be obtained following
  the proof of Theorem \ref{convrate1}.

  Recalling the
  proof of Theorem \ref{3dlom}, we actually showed that
  \be \|\mathbf{d}(t)-\mathbf{d}^*\|_{\mathbf{H}^1}\leq
  \varpi, \quad \forall\,
  t\geq 0,\label{esdd}
  \ee
   which implies that (let $t$ be large)
  $$\|\mathbf{d}_\infty-\mathbf{d}^*\|_{\mathbf{H}^1}\leq
  \|\mathbf{d}(t)-\mathbf{d}_\infty\|_{\mathbf{H}^1}+
  \|\mathbf{d}(t)-\mathbf{d}^*\|_{\mathbf{H}^1}<2\varpi\leq \min\{\beta, \delta\}.$$
  Taking $\mathbf{d}=\mathbf{d}_\infty$ and
  $\psi=\mathbf{d}^*$ in Corollary \ref{ELS},  we see that
  $$|\mathscr{E}(\mathbf{d}_\infty)-\mathscr{E}(\mathbf{d}^*)|^{1-\theta
  }=|\wha{E}(\mathbf{d}_\infty)-\wha{E}(\mathbf{d}^*)|^{1-\theta}\leq
  \|-\Delta \wha{\mathbf{d}}^*+\mathbf{f}(\mathbf{d}^*)\|=0.$$
  Since $\|\mathbf{d}_\infty-\mathbf{d}^*\|_{\mathbf{H}^1}\leq\delta$, $\mathbf{d}_\infty$ is also an energy minimizer of
  $\mathscr{E}$.

  Moreover, the proof of Theorem \ref{3dlom} implies that, for arbitrary (small) $\epsilon>0$, if we replace the choice of $\varpi$ \eqref{vpi} by
  \be
 \varpi_1=\min\left\{\epsilon, \frac{\beta}{2},\  \frac{\delta}{2},\
 \frac{\varepsilon_0}{10C_4}\right\},\label{vpi1}
 \ee
 then we are able to choose the constants $\sigma_i, M_j$ sufficiently small in a similar manner such that \eqref{conaaa} and \eqref{esdd} hold with $\varpi$ being replaced by $\varpi_1$ (and thus \eqref{esdd} holds for $\epsilon$). This yields the (locally) Lyapunov stability of $\mathbf{d}^*$. Finally, it is easy to see that if $\mathbf{d}^*$ is an isolated local
  minimizer, then  $\mathbf{d}_\infty=\mathbf{d}^*$ and $\mathbf{d}^*$ is asymptotically stable. The proof is complete.
  \end{proof}

\section{Appendix }
\label{App}
 \setcounter{equation}{0} We report
some properties of the lifting functions $\mathbf{d}_E$ and
$\mathbf{d}_P$ (cf. \eqref{LE} and \eqref{LP}) that have been
used in the previous sections. Below we denote by $c$ a generic
positive constant which depends on $n$ and $\Omega$ at most.

 \bl\label{Ap1}
 For any $t\geq 0$, and $k=0,1,2,...$, $j=0,1$, we have\\
 (i) $\|\p_t^j\mathbf{d}_E(t)\|_{\mathbf{H}^{k}}
 \leq c\|\partial_t^j\mathbf{h}(t)\|_{\mathbf{H}^{k-\frac12}(\Gamma)}
 $;\\
 (ii)
 $ \|\mathbf{d}_E(t)-\mathbf{d}_*\|_{\mathbf{H}^k}\leq c\|\mathbf{h}(t)-\mathbf{h}_*\|_{\mathbf{H}^{k-\frac12}(\Gamma)}$,
 where $\mathbf{d}_*$ is the unique solution to
 \be
 \begin{cases}
 -\Delta \mathbf{d}_*=\mathbf{0},\quad x\in \Omega,\\
 \mathbf{d}_*=\mathbf{h}_*, \quad x\in \Gamma.
 \end{cases}\label{dEI}
 \ee
 \el
 \begin{proof}
 The conclusion follows from the classical elliptic regularity theory (cf., e.g., \cite{LM,Tay}).
 \end{proof}


 \bl \label{Ap2}
 Let $\mathbf{d}_0\in \mathbf{H}^2(\Omega)$ with
 $|\mathbf{d}_0|_{\mathbb{R}^n}\leq 1$.
Suppose that $\mathbf{h}$ satisfy \eqref{hyp4}--\eqref{hyp5} and
$\mathbf{h}_t\in L^2_{loc} ([0,+\infty);
\mathbf{H}^{\frac12}(\Gamma))$.   Then, for any $t>0$, the following estimates hold
 \bea
  \| \mathbf{d}_P(t)-\mathbf{d}_E(t) \|^2_{\mathbf{H}^1}
 &\leq&  ce^{-t}\int_0^t e^\tau
 \|\mathbf{h}_t(\tau)\|_{\mathbf{H}^{-\frac12}(\Gamma)}^2d \tau ,
 \label{A3}
 \\
  \|\p_t\mathbf{d}_P(t)\|^2+ \|\mathbf{d}_P(t)-\mathbf{d}_E(t)\|_{\mathbf{H}^2}^2
 &\leq& c\int_0^t\|\mathbf{h}_t(\tau)\|^2_{\mathbf{H}^{\frac12}(\Gamma)}d\tau,\label{A5}\\
 \int_0^t \|\nabla \Delta \mathbf{d}_P\|^2 d\tau &\leq&
 c\int_0^t  \|\mathbf{h}_t(\tau)\|_{\mathbf{H}^{\frac12}(\Gamma)}^2
 d\tau.
 \label{A6}
 \eea
In addition, we have

\medskip\noindent
(i) if $\mathbf{h}_t\in L^2 (0,+\infty;
\mathbf{H}^{\frac12}(\Gamma))$ then
 \be
 \lim_{t\to+\infty}\|\partial_t \mathbf{d}_P(t)\|=0,\label{dedpt}
 \ee

\medskip\noindent
 (ii) if $\mathbf{h}_t$ satisfies (H6) then, for all $t\geq 0$,
 \bea
  \|\p_t\mathbf{d}_P(t)\|^2+
 \|\mathbf{d}_P(t)-\mathbf{d}_E(t)\|_{\mathbf{H}^2}^2&\leq& c
 (1+t)^{-2-2\gamma},  \label{A8}\\
  \int_{\frac{t}{2}}^t \|\nabla \Delta \mathbf{d}_P(\tau)\|^2 d\tau &\leq& c
 (1+t)^{-1-2\gamma}.\label{A9}
 \eea

 \el
 \begin{proof}
 It follows from  \eqref{LE} and \eqref{LP} that
 \be
 \begin{cases}
 -\Delta (\mathbf{d}_P-\mathbf{d}_E)=-\p_t\mathbf{d}_P,\qquad
 \text{ in }\Omega\times \mathbb{R}^+,\\
 \mathbf{d}_P-\mathbf{d}_E=\mathbf{0},\qquad \text{ on }
 \Gamma\times \mathbb{R}^+,
 \end{cases}\label{dEP}
 \ee
 and
  \be
 \begin{cases}
 \p_{t}(\mathbf{d}_P-\mathbf{d}_E)-\Delta
 (\mathbf{d}_P-\mathbf{d}_E)=-\p_{t}\mathbf{d}_E,\qquad \text{ in } \Omega\times \mathbb{R}^+,\\
 \mathbf{d}_P-\mathbf{d}_E=\mathbf{0},\qquad \text{ on } \Gamma\times \mathbb{R}^+,\\
 \mathbf{d}_P-\mathbf{d}_E|_{t=0}= \mathbf{0}, \qquad \text{ in }
 \Omega.
 \end{cases}
 \label{diEF}
 \ee
 Multiplying the first equation in \eqref{diEF} by $(\mathbf{d}_P-\mathbf{d}_E)-\Delta (\mathbf{d}_P-\mathbf{d}_E)$,
 integrating by parts and using the Poincar\'e inequality, we obtain
 \bea
 &&  \frac12\frac{d}{dt} (\| \mathbf{d}_P-\mathbf{d}_E \|^2+\|\nabla (\mathbf{d}_P-\mathbf{d}_E) \|^2)
  +\|\nabla (\mathbf{d}_P-\mathbf{d}_E) \|^2+\|\Delta (\mathbf{d}_P-\mathbf{d}_E) \|^2 \non\\
 &\leq& (\|\mathbf{d}_P-\mathbf{d}_E\| +\|\Delta (\mathbf{d}_P-\mathbf{d}_E) \|)\|\p_{t}\mathbf{d}_E\|\non\\
 &\leq& (C_P \|\nabla (\mathbf{d}_P-\mathbf{d}_E)\|+\|\Delta (\mathbf{d}_P-\mathbf{d}_E) \|) \|\p_{t}\mathbf{d}_E\|\non\\
 &\leq& \frac12 \|\nabla (\mathbf{d}_P-\mathbf{d}_E) \|^2+ \frac12 \|\Delta (\mathbf{d}_P-\mathbf{d}_E) \|^2
 + \left(\frac12 C_P^2
 +\frac12\right)\|\p_{t}\mathbf{d}_E\|^2,\label{dPE}
 \eea
 which, together with Lemma \ref{Ap1}, implies
 \bea
 \| \mathbf{d}_P(t)-\mathbf{d}_E(t) \|^2_{\mathbf{H}^1}
 &\leq&  ce^{-c_1t}\int_0^t e^{c_1\tau}\|\p_{t}\mathbf{d}_E(\tau)\|^2d\tau\non\\
  &\leq & ce^{-c_1t}\int_0^t e^{c_1\tau}
  \|\mathbf{h}_t(\tau)\|_{\mathbf{H}^{-\frac12}(\Gamma)}^2d \tau,
  \label{A3a}
 \eea
that is, \eqref{A3}.

Applying now the Laplacian to the first equation in \eqref{diEF}, we
get
  \be
 \begin{cases}
 \p_{t}\Delta(\mathbf{d}_P-\mathbf{d}_E)-\Delta^2 (\mathbf{d}_P
 -\mathbf{d}_E)=\mathbf{0},\qquad \text{ in }\Omega\times \mathbb{R}^+,\\
 \Delta(\mathbf{d}_P-\mathbf{d}_E)=\mathbf{h}_t,\qquad \text{ on }\Gamma\times \mathbb{R}^+,\\
 \Delta(\mathbf{d}_P-\mathbf{d}_E)|_{t=0}= \mathbf{0}, \qquad \text{ in }
 \Omega.\label{AAep}
 \end{cases}
 \ee
 Multiplying the first equation of \eqref{AAep}
 by $\Delta(\mathbf{d}_P-\mathbf{d}_E)$ and integrating by parts, we get
 \bea
  && \frac12\frac{d}{dt} \|\Delta(\mathbf{d}_P-\mathbf{d}_E)\|^2+ \|\nabla
  \Delta(\mathbf{d}_P-\mathbf{d}_E)\|^2\non\\
 &\leq&
 \|\partial_{\mathbf{n}}\Delta(\mathbf{d}_P-\mathbf{d}_E)\|_{\mathbf{H}^{-\frac12}(\Gamma)}
 \|\mathbf{h}_t\|_{\mathbf{H}^{\frac12}(\Gamma)}\non\\
 & \leq&  c\|\Delta (\mathbf{d}_P-\mathbf{d}_E)
 \|_{\mathbf{H}^1}\|\mathbf{h}_t\|_{\mathbf{H}^{\frac12}(\Gamma)}\non\\
 &\leq& \frac12 (\|\nabla
 \Delta(\mathbf{d}_P-\mathbf{d}_E)\|^2+\|\Delta(\mathbf{d}_P-\mathbf{d}_E)\|^2)+c
 \|\mathbf{h}_t\|_{\mathbf{H}^{\frac12}(\Gamma)}^2.\label{dPEt}
 \eea
 Hence, from \eqref{dPE} and \eqref{dPEt} we infer
  \be
\frac{d}{dt} \|\mathbf{d}_P(t)-\mathbf{d}_E(t)\|_{\mathbf{H}^2}^2+
c_2( \|\mathbf{d}_P(t)-\mathbf{d}_E(t)\|_{\mathbf{H}^2}^2+\|\nabla
\Delta (\mathbf{d}_P-\mathbf{d}_E)\|^2)\leq c
\|\mathbf{h}_t\|_{\mathbf{H}^{\frac12}(\Gamma)}^2,\label{dPE2}
  \ee
which entails \eqref{A6} and
   \be
\|\mathbf{d}_P(t)-\mathbf{d}_E(t)\|_{\mathbf{H}^2}^2
 \leq c\int_0^t
 \|\mathbf{h}_t(\tau)\|_{\mathbf{H}^{\frac12}(\Gamma)}^2
 d\tau.\label{A4}
\ee Thus \eqref{A5} follows from \eqref{A4} and the fact $
 \|\p_t\mathbf{d}_P(t)\| = \|\Delta \mathbf{d}_P(t)\|$.

Now if $\mathbf{h}_t\in L^2 (0,+\infty;
\mathbf{H}^{\frac12}(\Gamma))$, we infer from \eqref{dPE} that
 \bea
 \int_0^{+\infty} \|\Delta (\mathbf{d}_P(t)-\mathbf{d}_E(t))\|^2 dt
 &\leq&
 c\int_0^{+\infty}
 \|\p_{t}\mathbf{d}_E(t)\|^2dt\non\\
 &\leq& c\int_0^{+\infty}
  \|\mathbf{h}_t(t)\|_{\mathbf{H}^{-\frac12}(\Gamma)}^2
  dt<+\infty.\label{intdpt}
 \eea
  Then it follows from \eqref{dPEt}, \eqref{intdpt} and Lemma
  \ref{SZ} that
  $$
  \lim_{t\to +\infty} \|\Delta
  (\mathbf{d}_P(t)-\mathbf{d}_E(t))\|^2=0,$$  which implies
  \eqref{dedpt}.

 Furthermore, if  (H6) holds, then \eqref{dPE2} implies that (cf., e.g., \cite{WGZ1})
    \be
   \|\mathbf{d}_P(t)-\mathbf{d}_E(t)\|_{\mathbf{H}^2}^2\leq
   c(1+t)^{-2-2\gamma}, \quad \forall \ t\geq 0.\non
    \ee
    Using \eqref{dPE2} once more, we have
    \bea
    &&\int_{\frac{t}{2}}^t \|\nabla \Delta \mathbf{d}_P(\tau)\|^2 d\tau =
    \int_{\frac{t}{2}}^t \|\nabla \Delta (\mathbf{d}_P-\mathbf{d}_E)(\tau)\|^2 d\tau\non\\
    & \leq&
    c \left\|\mathbf{d}_P\left(\frac{t}{2}\right)-\mathbf{d}_E\left(\frac{t}{2}\right)\right\|_{\mathbf{H}^2}^2
    + c\int_{\frac{t}{2}}^t  \|\mathbf{h}_t(\tau)\|_{\mathbf{H}^{\frac12}(\Gamma)}^2
    d\tau\non\\
    &\leq& c\left(1+\frac{t}{2}\right)^{-2-2\gamma}+\frac{c}{1+2\gamma}\left(1+\frac{t}{2}\right)^{-1-2\gamma}\non\\
    & \leq& c\left(1+t\right)^{-1-2\gamma}, \quad \forall\ t\geq 0,\non
    \eea
    and this gives \eqref{A9}. The proof is complete.
  \end{proof}

\medskip

\noindent {\bf Acknowledgments.} The authors would like to thank the referees
for their helpful comments and suggestions on an earlier
version of this paper. This work originated from a visit
of the first author to the Fudan University whose hospitality is
gratefully acknowledged. The second author
was partially supported by NSF of China 11001058, Specialized Research Fund for
the Doctoral Program of Higher Education and "Chen Guang" project supported by Shanghai Municipal Education Commission and Shanghai Education Development Foundation.


\end{document}